%% file: main.tex
\documentclass[final]{siamltex}

\usepackage{graphicx}
\usepackage{amssymb}
\usepackage{amsmath}
\usepackage{subfig}
\usepackage{multirow}
\usepackage{colortbl}
\usepackage{algorithm}
\usepackage{algpseudocode}

\providecommand{\slfrac}[2]{\left.#1\middle/#2\right.}

\title{The cost of continuity: performance of iterative solvers on isogeometric finite elements}

\author{Nathan~Collier\thanks{Center for Numerical Porous Media
    (NumPor), Applied Mathematics and Computational Science, Earth and
    Environmental Sciences and Engineering, King Abdullah University
    of Science and Technology, Thuwal, Saudi Arabia ({\tt
      nathaniel.collier@gmail.com})} \and
  Lisandro~Dalcin\thanks{Consejo Nacional de Investigaciones
    Cient\'\i{}ficas y T\'ecnicas, Santa Fe, Argentina. This author
    was partially funded by the Center for Numerical Porous Media,
    KAUST Strategic Research Initiative} \and
  David~Pardo\thanks{Department of Applied Mathematics, Statistics,
    and Operational Research, The University of the Basque Country
    UPV/EHU and Ikerbasque, Bilbao, Spain. This author was partially
    funded by the Project of the Spanish Ministry of Sciences and
    Innovation MTM2010-16511, the Laboratory of Mathematics (UFI
    11/52)} \and V.~M.~Calo\thanks{Co-director, Center for Numerical
    Porous Media (NumPor), Applied Mathematics and Computational
    Science, Earth and Environmental Sciences and Engineering, King
    Abdullah University of Science and Technology, Thuwal, Saudi
    Arabia}}

\begin{document}

\maketitle

\begin{abstract}
  In this paper we study how the use of a more continuous set of basis
  functions affects the cost of solving systems of linear equations
  resulting from a discretized Galerkin weak form. Specifically, we
  compare performance of linear solvers when discretizing using $C^0$
  B-splines, which span traditional finite element spaces, and
  $C^{p-1}$ B-splines, which represent maximum continuity. We provide
  theoretical estimates for the increase in cost of the matrix-vector
  product as well as for the construction and application of black-box
  preconditioners. We accompany these estimates with numerical results
  and study their sensitivity to various grid parameters such as
  element size $h$ and polynomial order of approximation $p$. Finally,
  we present timing results for a range of preconditioning options for
  the Laplace problem. We conclude that the matrix-vector product
  operation is at most $\slfrac{33p^2}{8}$ times more expensive for
  the more continuous space, although for moderately low $p$, this
  number is significantly reduced. Moreover, if static condensation is
  not employed, this number further reduces to at most a value of 8,
  even for high $p$. Preconditioning options can be up to $p^3$ times
  more expensive to setup, although this difference significantly
  decreases for some popular preconditioners such as Incomplete LU
  factorization.
\end{abstract}


\begin{AMS}\end{AMS}

\pagestyle{myheadings}
\thispagestyle{plain}
\markboth{COLLIER, DALCIN, PARDO AND CALO}{IGA \& ITERATIVE SOLVERS}

\section{Introduction}\label{s:intro}\input{intro.tex}
\section{Model Problem}\label{s:model_problem}\input{model.tex}
\section{Theory}\label{s:theory}\input{theory.tex}
\section{Numerical Results}\label{s:res}\input{res.tex}
\section{Conclusions}\label{s:conc}\input{conc.tex}

\section*{Acknowledgments}
We would like to thank Analissa Buffa and Juan Galvis for the
interesting discussions we had with them on this topic. Also, we would
like to acknowledge the open source software packages that made this
work possible: {\tt PETSc}~\cite{petsc1,petsc2}, {\tt
  petsc4py}~\cite{Dalcin2011}, {\tt NumPy}~\cite{numpy}, {\tt
  matplotlib}~\cite{matplotlib}, {\tt SymPy}~\cite{sympy}. 

\bibliographystyle{siam} \bibliography{paper}

\end{document}

%% file: intro.tex
Isogeometric analysis (IGA)~\cite{Hughes2005,Cottrell2009} is a
Galerkin finite element method which has popularized the use of the
Non-uniform Rational B-spline (NURBS) basis for solving partial
differential equations (PDE's). The computer aided design (CAD)
community has long used NURBS as a basis due to its higher-order
continuity, ideal for designing smooth curves and surfaces. The higher
order continuous basis also enables the use of the standard Galerkin
method for solving higher order
problems~\cite{Gomez2008,Gomez2010,Dede2010}. Furthermore, it has been
observed that the approximability of the higher continuous spaces per
degree of freedom is superior to that of traditional finite element
spaces for problems with smooth solutions. This suggests that IGA not
only links geometry to analysis but also is an efficient method for
solving a variety of PDE's.

As for any Galerkin method, the main computational cost of IGA comes
from the assembly and solution of a system of linear equations. When
using a direct method to solve this linear system, we showed
in~\cite{Collier2012} that for large three dimensional problems of
given polynomial order $p$, the number of floating point operations
(FLOPS) needed to solve a $C^{p-1}$ discretization is approximately
$p^3$ times more expensive than that of a $C^0$ discretization with
the same order of approximation $p$ and same number of degrees of
freedom (DOF). In other words, each DOF is $p^3$ times more expensive
to solve when using $C^{p-1}$ discretizations as opposed to using
$C^{0}$ discretizations. This theoretical estimate was corroborated by
our numerical experiments for problem sizes of practical interest.

In this work, we extend our previous study to the case of iterative
solvers. The main motivation for using iterative methods is to reduce
computational cost (time and memory). However, in general there are no
{\em a priori} estimates for the computational time required by the
iterative solver because it is a composition of many factors. These
costs include matrix-vector multiplications and additional operations
required by the iterative method (vector scalings, dot products) as
well as the cost of setting up and applying the preconditioner. While
these costs may be estimated, their influence on the overall
computational time tightly depends on the number of iterations
required to reduce the linear algebra error to a prescribed
tolerance. For example, it often occurs that a preconditioner which
leads to few required iterations for convergence is also more
expensive to construct and/or apply. In the limit, a LU factorization
is the ideal preconditioner in terms of iteration count, however is
expensive to construct and apply.

To develop a baseline understanding for how continuity affects
iterative solvers, we study the canonical Laplace problem discretized
using $C^0$ and $C^{p-1}$ B-spline spaces, representing minimum and
maximum continuity. We only consider the matrix-vector multiplication
component of the iterative solver. The additional operations required
by different iterative methods (vector updates, orthogonalizations)
are not dependent on the basis used, and therefore may be ignored when
comparing how continuity affects the iterative solver. To better
expose the effect of continuity on the cost of the solver, we use
preconditioned conjugate gradients (CG) as our iterative method,
because it is among the most efficient methods in the Krylov
family~\cite{Saad,CGnopain}.

We will study a range of standard preconditioners which are
appropriate for small and medium size problems. These include:
diagonal Jacobi, successive overrelaxation, incomplete LU
factorization, and element by element. Despite the fact these methods
do not scale to large problems, they are frequently used as building
blocks to construct more complex preconditioners (approximate solvers
on subdomains of domain decomposition and physics-based
preconditioners or smoothers in multigrid techniques
\cite{Pardo2004,Ainsworth1996,Arnold2000}). Most of the techniques we
will study here are described in Saad's book \cite{Saad} and
implemented in scientific software frameworks such as
PETSc~\cite{petsc1,petsc2}. It is our aim to assess the additional
cost incurred by a more continuous basis as well as illuminate how
standard approaches work for IGA discretizations.

This study complements the only previous published work on iterative
solvers for IGA discretizations of which we are aware. In
\cite{Veiga2012tr} Beir\~{a}o~da~Veiga, et. al. propose a large family
of domain-decomposition two-grid solvers and prove theoretically that
the condition number of the preconditioned system is independent of
element size $h$. They also provide numerical evidence showing that
the condition number is independent of $p$ provided that the overlap
between subdomains is sufficiently large. However, these numerical
results are concerned only with convergence (condition number and
number of iterations) and not with computational efficiency.

The rest of this paper develops in the following manner. In
section~\ref{s:model_problem} we detail the model problem used
throughout this work. Section~\ref{s:theory} derives theoretical
estimates of FLOPS needed to perform matrix-vector multiplications of
linear systems resulting from $C^0$ and $C^{p-1}$ discretizations, as
well as estimates for the setup and application of different
preconditioners. In section~\ref{s:res} we present numerical results
to complement the theory. We show convergence in terms of iterations
as well as computational time on a range of discretizations varying in
$h$ and $p$.

%% file: model.tex
The problem used for our study is the Laplace equation in three
dimensions on the unit cube,
\begin{equation}
  \begin{aligned}
    -\nabla\cdot(\nabla u)&=0\ \ \ \text{on }\Omega\\
    u&=1\ \ \ \text{on }\Gamma_D\\
    (\nabla u)\cdot\mathbf{n}&=0\ \ \ \text{on }\Gamma_N
  \end{aligned}
\end{equation}
where $\Omega=[0,1]^3$, $\Gamma_D=(0,:,:) \cup (:,0,:) \cup (:,:,0)$,
and $\Gamma_N=(1,:,:) \cup (:,1,:) \cup (:,:,1)$. We will use uniform
$h$-refinements of $C^0$ and $C^{p-1}$ B-splines to discretize the
weak form of the Laplace equation.

%% file: theory.tex
In this section we develop theoretical estimates for the increase in
cost associated with the use of higher continuous spaces in Galerkin
finite elements. We assess cost by counting the FLOPS required by
matrix-vector products and the setup of different preconditioning
options. We use these estimates as a measure of the relative cost
between $C^0$ and $C^{p-1}$ spaces.

\subsection{Matrix-Vector Multiplication}

The main cost of iterative methods is due to the matrix-vector
multiplication operation which is proportional to the number of
nonzero entries in both the system and preconditioner matrix. We
develop estimates for the number of nonzero entries in the stiffness
matrix of the model problem resulting from $C^0$ and $C^{p-1}$
discretizations in three spatial dimensions. We do this by considering
the number of nonzero entries that a single element of a structured
grid mesh contributes to the system matrix. 

\begin{figure}[ht]
\centering
\includegraphics[width=0.6\textwidth]{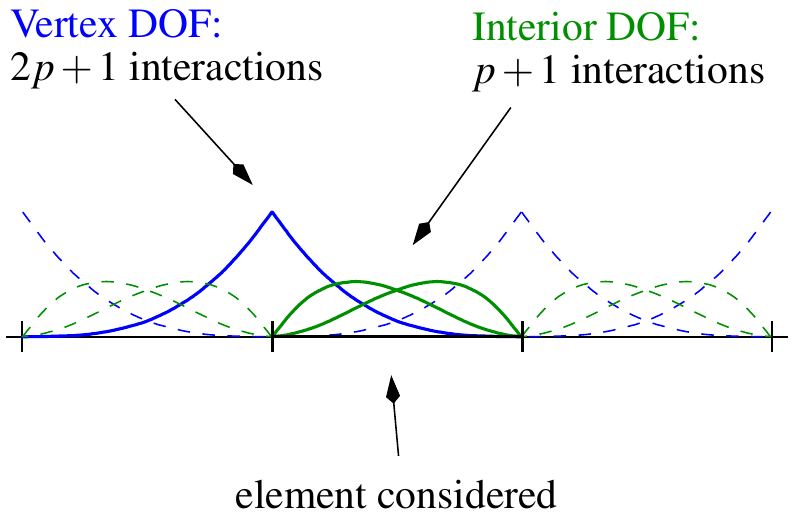}
\caption{Sample cubic $C^0$ discretization with a single element
  highlighted.}\label{f:nnz1d}
\end{figure}

We begin by considering a single element of a 1D, $C^0$ discretization
of order $p$. Consider figure~\ref{f:nnz1d} where we have drawn such
an element, particularized to a cubic for the sake of
illustration. The number of nonzero entries to the matrix will be the
sum of the interactions that each basis has with all other basis
functions which have overlapping support. We note that in the 1D case,
there are two classes of interactions--those associated to the
vertices and the interior of the element. We note that a basis
associated to a vertex overlaps $2p+1$ others while the bases
associated to interiors overlap $p+1$ others. We consider a DOF for a
single vertex per element (since the other vertex is actually the
first vertex of the next element) and $p-1$ DOFs per interior. The
total number of nonzero entries accumulated due to a single element is
then $(1)(2p+1)+(p-1)(p+1)$. This number is attained by summing over
all entities the total number of interactions of all the DOFs
associated to that entity.

For the multidimensional case, we extend this method of counting the
nonzero entries of the matrix by a tensor product construction. In
addition to vertex and interior DOFs, we have DOFs associated to edges
in two and three dimensions as well as DOFs associated to faces in
three dimensions. We summarize the enumeration of DOFs and their
interactions in table~\ref{t:interact}. The total number of nonzeros
due to a single element contribution is the row sum of the product of
all columns in this table. Specially, for three dimensions we have
\begin{equation}
\begin{tabular}{rcl}
$\mbox{nnz}^{C^0} $ & $=$ & $\underbrace{(p-1)^3}_{\mbox{interior DOF}} \cdot\  (p+1)^3$\\
[.35in]
 & $+$ & $\underbrace{3(p-1)^2}_{\mbox{face DOF}} \cdot\  (2p+1) (p+1)^2$  \\
[.35in]
 & $+$ & $\underbrace{3 (p-1)}_{\mbox{edge DOF}} \cdot\  (2p+1)^2(p+1)$ \\
[.35in]
 & $+$ & $\underbrace{1}_{\mbox{vertex DOF}}\cdot\  (2p+1)^3$\\
[.35in]
 & $=$ & $p^6 + 6p^5 + 12p^4+8p^3$\\ 
[.35in]
 & $=$ & $p^3(p+2)^3 = {\cal O}(p^6)$ 
\end{tabular}\label{e:nnzc0}
\end{equation}
\begin{table}
\centering
\caption{Summary table of the interactions of degrees of freedom
  associated with a $C^0$ basis.}\label{t:interact}
\begin{tabular}{lllll}
\hline
 &  & Number & DOFs & Number  \\
Dimension & Entity & of Entities & per Entity & of interactions\\
\hline
1D & vertex & 1 & 1 & $(2p+1)$\\
1D & interior & 1 & $(p-1)$ & $(p+1)$\\
\rowcolor[gray]{.9}2D & vertex & 1 & 1 & $(2p+1)^2$\\
\rowcolor[gray]{.9}2D & edge & 2 & $(p-1)$ & $(2p+1)(p+1)$\\
\rowcolor[gray]{.9}2D & interior & 1 & $(p-1)^2$ & $(p+1)^2$\\
3D & vertex & 1 & 1 & $(2p+1)^3$\\
3D & edge & 3 & $(p-1)$ & $(2p+1)^2(p+1)$\\
3D & face & 3 & $(p-1)^2$ & $(2p+1)(p+1)^2$\\
3D & interior & 1 & $(p-1)^3$ & $(p+1)^3$\\
\hline
\end{tabular}
\end{table}

In the case of $C^{p-1}$ B-splines, the interactions are more regular
because each basis interacts with $(2p+1)^3$ others. To make the
estimates comparable, in terms of unknowns, to the single element of
$C^0$ basis functions, we multiply by $p^3$.
\begin{equation}
\mbox{nnz}^{C^{p-1}} = p^3(2p+1)^3 = 8p^6 + 12p^5 + 6p^4 + p^3 = {\cal O}(8p^6)\label{e:nnzcpm1}
\end{equation}

We conclude that in the case of large $p$ the increase in cost of
matrix vector multiplication of $C^{p-1}$ spaces is no more than eight
times that of $C^0$ spaces. However, for the range of meaningful
discretizations of polynomial order $p$, we see this factor smaller
than the limit, approximately two for $p=2$ and three for $p=3$. In
table~\ref{tab:memdiff} we present some numerical results for the
actual ratios of times for 1000 matrix-vector products of $C^{p-1}$
and $C^0$ spaces as the number of degrees of freedom, $N$, in the
system increases. We compare these time ratios to the theoretical
ratio of number of nonzero entries for a system of infinite size,
equation~\eqref{e:nnzcpm1} divided by equation~\eqref{e:nnzc0}.

\begin{table}[htp]
\centering
\caption{Actual and estimated increase in cost of a matrix-vector
  multiply for $C^{p-1}$ spaces relative to $C^0$ spaces.}
\label{tab:memdiff}
\begin{tabular}{ccccc}
\hline
 & \multicolumn{4}{c}{Polynomial order, $p$}\\
\cline{2-5}
$N$ & $2$ & $3$ & $4$ & $5$ \\
\hline
$10^3$ & 1.19 & 1.95 & 1.93 & 0.94 \\
$10^4$ & 1.76 & 2.22 & 2.96 & 3.19 \\
$10^5$ & 1.74 & 2.56 & 3.02 & 3.46 \\
$10^6$ & 1.80 & 2.60 & 3.08 & 3.51 \\
$\infty$ & 1.95 & 2.74 & 3.37 & 3.88 \\
\hline
\end{tabular}
\end{table}

\subsubsection*{Note on Static Condensation}

When using $C^0$ spaces, it is common to first eliminate (using
Gaussian elimination) all DOF interior to an element, a technique
known as {\em static condensation}~\cite{Wilson1974}. This approach is
also used in a multi-frontal direct solver algorithm
\cite{Duff1983,Duff1984} and known to be of reduced value when using
$C^{p-1}$ spaces (see~\cite{Collier2012}). Iterative solvers can also
make use of the technique, solving on the reduced system, called the
skeleton problem. The skeleton problem is not only of smaller rank
than the original, but it also contains fewer nonzero entries. This
affects the iterative solver in that the matrix-vector multiplications
are economized.

To see this effect, we compute the number of nonzeros for a single
element in the resulting matrix after performing static
condensation. We repeat a portion of table~\ref{t:interact} which
corresponds to the three dimensional results in
table~\ref{t:intsc}. If we statically condense the interior DOFs,
these nonzero entries are now removed (the row of the table is
removed). However, we also need to remove all interactions that the
vertices, edges, and faces have with interior DOFs. To this end, we
have added another column which represents these DOFs. For each entity
we eliminate $(p-1)^3$ DOFs which correspond to each interior to which
that entity was connected. Vertices connect to eight interiors, edges
to four interiors, and faces to two interiors. We then sum the nonzero
entries as before.
\begin{table}
\centering
\caption{Summary table of the interactions of degrees of freedom
  associated with a $C^0$ basis in three dimensions with the interior
  DOFs statically condensed.}\label{t:intsc}
\begin{tabular}{lllll}
\hline
 &  Number & DOFs & Number & Statically \\
Entity & of Entities & per Entity & of interactions & condensed\\
\hline
vertex & 1 & 1 & $(2p+1)^3$ & $-8(p-1)^3$ \\
edge & 3 & $(p-1)$ & $(2p+1)^2(p+1)$ & $-4(p-1)^3$ \\
face & 3 & $(p-1)^2$ & $(2p+1)(p+1)^2$ & $-2(p-1)^3$ \\
\hline
\end{tabular}
\end{table}
\begin{equation}
\begin{tabular}{rcl}
$\mbox{nnz}_{sc}^{C^0}$ & $=$ & $\underbrace{3(p-1)^2}_{\mbox{face DOF}} \cdot\  [(2p+1)(p+1)^2-2(p-1)^3]$  \\
[.35in]
 & $+$ & $\underbrace{3 (p-1)}_{\mbox{edge DOF}} \cdot\  [(2p+1)^2(p+1)-4(p-1)^3] $ \\
[.35in]
 & $+$ & $\underbrace{1}_{\mbox{vertex DOF}} \cdot\  [(2p+1)^3-8(p-1)^3]$\\
[.35in]
 & $=$ & $ 33p^4 - 12 p^3 + 9p^2 -6p +3 = {\cal O}(33p^4)$ 
\end{tabular}
\end{equation}

For large enough problems, the matrix-vector product of $C^0$ spaces
becomes $\slfrac{p^2}{33}$ more expensive that the statically
condensed system. In the case of $C^{p-1}$ spaces, this number
approaches $\slfrac{8p^2}{33}$. While the process of static
condensation incurs additional cost in the matrix assembly phase, in
practice this approach is more efficient than standard $C^0$
approaches and not worthwhile in the case of $C^{p-1}$ spaces. See
table~\ref{t:memdiffsc} for comparison of theory to timing results.

\begin{table}[htp]
\centering
\caption{Actual and estimated increase in cost of a matrix-vector
  multiply for $C^{p-1}$ spaces relative to $C^0$ spaces with static
  condensation.}
\label{t:memdiffsc}
\begin{tabular}{ccccc}
\hline
 & \multicolumn{4}{c}{Polynomial order, $p$}\\
\cline{2-5}
$N$ & $2$ & $3$ & $4$ & $5$ \\
\hline
$10^3$ & 1.32 & 2.66 & 3.26 & 2.00 \\
$10^4$ & 1.96 & 3.07 & 5.14 & 6.87 \\
$10^5$ & 1.95 & 3.56 & 5.29 & 7.58 \\
$10^6$ & 2.00 & 3.62 & 5.42 & 7.74 \\
$\infty$ & 2.16 & 3.81 & 5.93 & 8.54 \\
\hline
\end{tabular}
\end{table}

\begin{figure}[ht]
\centering
\includegraphics[width=0.8\textwidth]{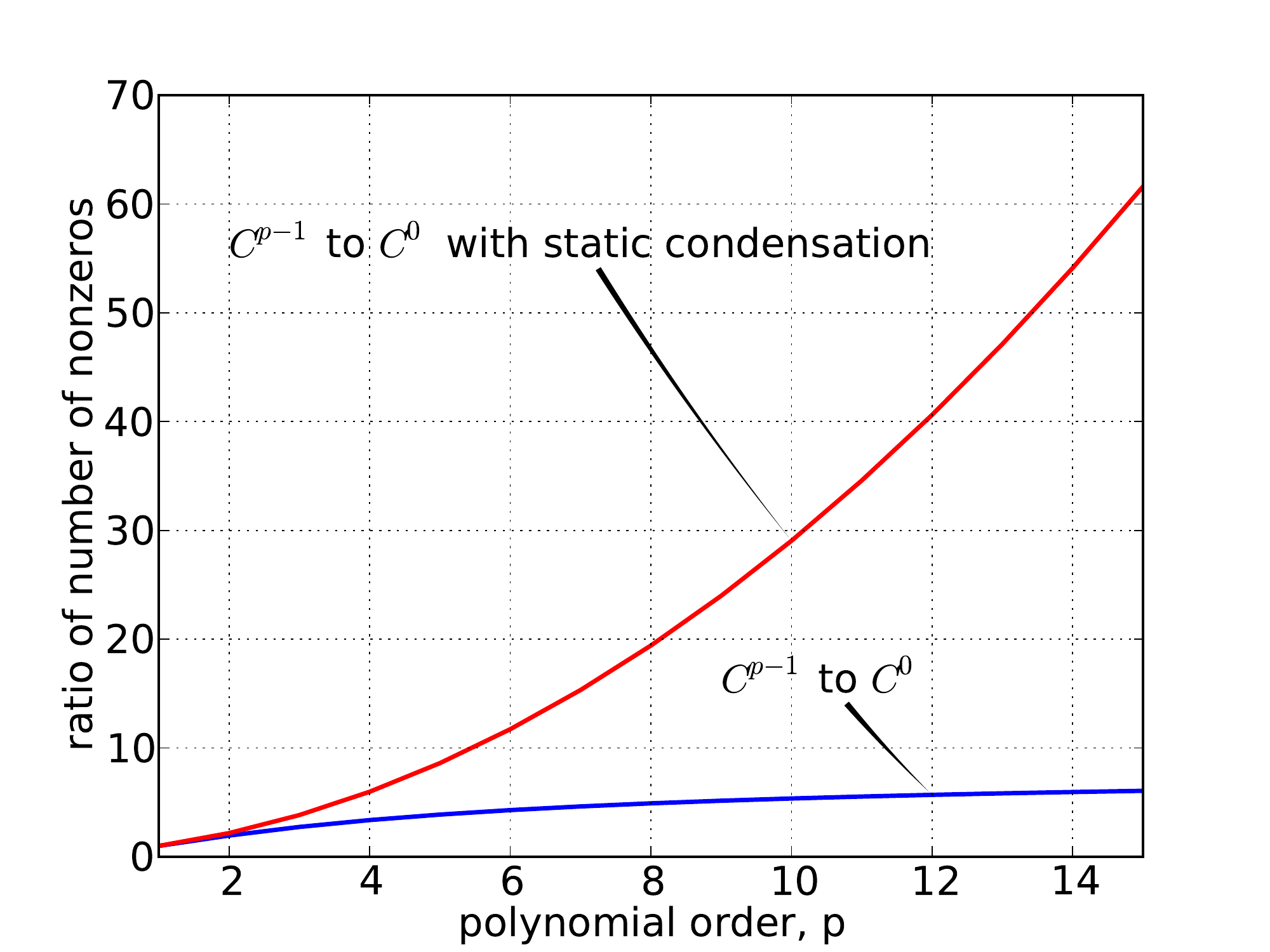}
\caption{The theoretical ratio of the number of nonzeros in system
  matrices, independent of the number of unknowns.}\label{f:nnz}
\end{figure}

While the gains in static condensation when using relatively low $p$
are moderate, for higher $p$ the added efficiency is of greater
importance. In figure~\ref{f:nnz} we plot the theoretical ratios of
the number of nonzero entries for $C^{p-1}$ relative to $C^0$ spaces
with and without static condensation. These plots represent, as $p$
increases, how much more expensive a matrix-vector product is when
using $C^{p-1}$ spaces. When no static condensation is used, we see
that the increase asymptotically approaches (slowly) a factor of
eight. However, when compared to the use of static condensation, the
increase in cost continues to grow with high $p$. If one is to
advocate the use of $C^{p-1}$ basis functions as a high $p$ method,
the merits of the basis must be weighed against this increase in cost.

\subsection{Black-box Preconditioners}

Now we consider how more continuous bases affect black-box
preconditioning techniques, such as those found in \cite{Saad}
(Chapter 10, pages 317--339). We develop theoretical cost estimates in
terms of FLOPS for both forming and applying each preconditioner. In
the paragraphs to follow we briefly describe each preconditioner and
explain how these estimates may be formed.
\subsubsection*{Diagonal Jacobi}
Practical implementations of diagonal Jacobi preconditioning extract
the diagonal entries from the matrix and invert them, storing the
result in a vector. The application of the preconditioner is then
performed by point-wise multiplication of residual entries with the
diagonal inverses. Both the setup and application of this
preconditioner require $N$ FLOPS, independent of the continuity of the
basis.
\subsubsection*{Symmetric Successive Overrelaxation (SSOR)}
The SSOR preconditioner is based on a relaxation scheme, similar to
Gauss-Seidel iterations. Practical implementations of SSOR
preconditioning extract the diagonal entries from the matrix, invert,
and scale them by the relaxation parameter in order to make the
application of the preconditioner more economical. The application of
this preconditioner consists of forward and backward sweeps, which
roughly amounts to a single matrix-vector product.

\subsubsection*{Incomplete LU factorization (ILU)}
Incomplete LU factorization (ILU) is a preconditioning technique based
on Gaussian elimination. Here we address only ILU with zero
fill-in. The ILU preconditioner is formed by performing LU
factorization, omitting entries which would change the nonzero pattern
of the original matrix. Thus, ILU is a crude approximation to the LU
factors of the system matrix, however more economical to compute and
apply.

Implementations of the zero fill-in ILU preconditioner are based in
the $IKJ$ (see the discussion in \cite{Saad} starting on page 304)
version of Gaussian elimination on the static non-zero pattern of the
input sparse matrix. The algorithm traverses the sparse matrix by
rows. At each row, the Gaussian elimination algorithm is applied on
only the nonzero entries.

Denoting $L_i$ the number of nonzero entries in the strictly
lower-triangular part of the $i$-th row and $U_k$ the number of
nonzero entries in the strictly upper-triangular part of the $k$-th
row, the number of FLOPs required to eliminate the $i$-th row is
$L_i(1+2U_k)$.

For a $C^{p-1}$ system matrix, every row has $(2p+1)^3$ nonzero
entries, thus the the number of nonzero entries in the strictly
lower-triangular and upper-triangular parts of $i$-th and $k$-th rows
are 
\[L_i = U_k = \frac{(2p+1)^3-1}{2}\]
and the total number of FLOPS for $N$ rows is
\[\text{FLOPS}^{C^{p-1}}_{\text{ILU}} = N \left(32 p^6 + 96 p^5 + 120 p^4 + 76 p^3 + 24 p^2 + 3 p \right)\]

For a $C^0$ system matrix, the number of nonzeros per row depends on
the kind of DOF (see previous subsection) and obtaining analytic
estimates is much more involved. We use instead a computational
approach consisting on building the graph for a mesh of $5 \times 5
\times 5$ elements for a $C^0$ discretization of degree
$p=1\ldots7$. For every $p$, we compute the preconditioner row-by-row
and add-up the number of FLOPS required for performing the ILU
factorization for the middle element. By using polynomial fitting, we
obtain the coefficients of a degree 6 polynomial. Finally, the cost of
constructing the ILU preconditioner for a $C^0$ system matrix is

\[\text{FLOPS}^{C^{0}}_{\text{ILU}} = N\left(\frac{2}{3} p^{6} + \frac{26}{3} p^{5} + \frac{128}{3} p^{4} + \frac{601}{6} p^{3} + \frac{355}{3} p^{2} + \frac{200}{3} p + \frac{83}{6}\right)\]

We highlight that the ILU preconditioner is inexpensive relative to
the full LU factorization and is in both cases of order $N$. We also
note that in three spatial dimensions, the highest order term in terms
of the polynomial order is $p^6$ for both $C^0$ and $C^{p-1}$
spaces. In the case of large $p$ the increase in cost of matrix vector
multiplication of $C^{p-1}$ spaces is no more than 48 times that of
$C^0$ spaces. In the summary table~\ref{t:cost}, we only include these
leading terms in order to more succinctly compare preconditioners. The
application of this preconditioner consists of forward and backward
substitution steps, which roughly amounts to a single matrix-vector
product.

\subsubsection*{Element-by-element (EBE)}

In the context of $C^0$ finite element spaces, an additive-Schwarz
\cite{Smith1996} preconditioner may be constructed in the limit where
the subdomains are the elements of the finite element
discretization. This preconditioner is known as the element-by-element
preconditioner (EBE). Note that this preconditioner departs from that
described in \cite{Saad} and follows \cite{Pardo2004}. Given the fully
assembled system matrix, this preconditioner is constructed by
extracting the local element matrix and inverting it explicitly. The
inverse of the local element matrix is assembled into a preconditioner
matrix which has the same nonzero pattern as the original system
matrix.

The cost of constructing the preconditioner for both $C^0$ and
$C^{p-1}$ spaces is the number of elements, $N_e$, times the cost of
inverting the small blocks,
\[N_e\left(2p^{9}\right)\] 
We note that for $C^0$ discretizations the number of degrees of
freedom $N$ can be related to the number of elements by the
relationship, $N=\mathcal{O}(N_ep^3)$. Therefore, the preconditioner
cost can then be expressed in terms of number of degrees of freedom
as, $2Np^{6}$. In $C^{p-1}$ spaces, the number of elements is roughly
the number of degrees of freedom, $N=\mathcal{O}(N_e)$ which leads to
the total cost being $2Np^{9}$. We emphasize that in this case, the
resulting matrix-vector product is no more expensive than for the
original system matrix. This means that the EBE preconditioner is
again at most 8 times more expensive to apply for $C^{p-1}$ spaces
when compared to $C^0$.

\subsubsection*{Basis by Basis (BBB)}

For $C^{p-1}$ spaces, we also construct an additive-Schwarz type
preconditioner based on the family of preconditioners presented in
\cite{Veiga2012tr}. We consider a selection of these preconditioners
constructed by taking single basis function subsets of the function
space. We explore the performance of this family of preconditioners by
varying the number of overlapping basis functions, $0\le r \le p$. We
call this preconditioner Basis by Basis (BBB) and note that if $r=0$,
the preconditioner corresponds to diagonal Jacobi. In the case that
$r=\slfrac{p}{2}$ the preconditioner is similar to the EBE
preconditioner described in this section. If the polynomial order is
even and the domain is periodic, it is identical to the element-based
preconditioner.

The cost of constructing this preconditioner is the number of basis
functions, multiplied by the cost of inverting the block,
$N2\left(2r+1\right)^{9}$. However, in this more general family of
preconditioners, the nonzero structure of the preconditioner matrix
varies with choice of $r$, resulting in a significant change in the
cost of the matrix-vector product. The number of nonzero entries in
the system matrix for a $C^{p-1}$ basis is $N(2p+1)^{9}$. The number
of nonzero entries in the preconditioner matrix can be obtained by a
similar expression, this time each row interacting with $4r+1$
columns. This leads to a number of nonzero entries which grows like
$N(4r+1)^{9}$. Thus the cost of the matrix-vector product of the
preconditioner matrix relative to the system matrix can be expressed
as $(\slfrac{2r}{p})^9$. If $r=\slfrac{p}{2}$, then applying the
preconditioner is just as expensive as a matrix-vector product of the
system matrix.

\subsubsection*{Summary}

We summarize the cost of setting up and applying each preconditioner
in table~\ref{t:cost}. We note that in all cases, except for the
trivial diagonal Jacobi or SSOR, the setup cost of these
preconditioners is more expensive for the $C^{p-1}$ spaces. Also, the
application of the preconditioners we study is in most cases no more
expensive than the matrix-vector product of the corresponding
space. Of particular interest in the case of $C^{p-1}$ spaces, is that
the EBE and BBB preconditioners are estimated to take $p^3$ more FLOPS
to setup, suggesting that they might not be as useful from a practical
point-of-view. This estimation is corroborated by the numerical
experiments in section~\ref{s:res}.

\begin{table}
\caption{Comparison of FLOPS estimates for the setup and application
  of different preconditioners for $C^0$ and $C^{p-1} $ spaces.}\label{t:cost}
\begin{center}
\begin{tabular}{llll}
\hline
Type\ \ \ \ \ & Space\ \ \ \ \ & Setup FLOPS\ \ & Apply FLOPS\ \ \ \\
\hline
Jacobi & $C^0$ & $N$ & $N$ \\
Jacobi & $C^{p-1}$ & $N$ & $N$ \\
\rowcolor[gray]{.9}SSOR & $C^0$ & $N$ & $2Np^3$\\
\rowcolor[gray]{.9}SSOR & $C^{p-1}$ & $N$ & $16Np^3$\\
ILU & $C^0$ & $\frac{2}{3}Np^6$ & $Np^3$\\
ILU & $C^{p-1}$ & $32Np^6$ & $8Np^3$\\
\rowcolor[gray]{.9}EBE & $C^0$ & $2Np^6$ & $Np^3$\\
\rowcolor[gray]{.9}EBE & $C^{p-1}$ & $2Np^9$ & $8Np^3$\\
BBB & $C^{p-1}$ & $2^{10}Nr^9$ & $(\slfrac{2r}{p})^98Np^3$\\
\hline
\end{tabular}
\end{center}
\end{table}

%% file: res.tex
In this section, we first present numerical results confirming the
theoretical estimates presented in the previous section. We do this to
isolate the cost of sparse matrix kernels from the iterative method in
which they are employed. Second, we present results on iteration
counts as the spaces are scaled in $h$ and $p$ for a range of
preconditioners. Finally, we report wall clock times required to solve
the model problem. In all our numerical tests, we start from an
initial guess of zero, and declare convergence when the preconditioned
residual norm decreases by eight orders of magnitude.

\subsection{Sparse Matrix Kernels}

In this subsection, we chose a periodic three-dimensional grid of $N =
60 \times60 \times 60 = 216,000$ degrees of freedom. By using
polynomial degrees $p = 1\dots 5$, we can construct $C^0$ and
$C^{p-1}$ discretizations with $(60/p) \times (60/p) \times (60/p)$
and $60 \times 60 \times 60$ elements respectively. For these
discretizations, we assemble consistent mass matrices and use them for
our experiments.

The experiments consist in measuring the wall-clock time spent in the
various sparse matrix kernels discussed in section~\ref{s:theory}. We
recall that sparse matrix-vector product operations, symmetric
successive over relaxations sweeps, and triangular solves require 2
FLOPS per nonzero entry in the sparse matrix.

The experiments were conducted on a desktop machine with a 3.07 GHz
Intel Core i7 950 processor, 8 megabytes of L3 cache and a cache line
of 64 bytes, a 4.8 GT/s Intel QPI memory interconnect, and 12
gigabytes of 1066 MHz DDR3 memory. The standard STREAM Triad benchmark
performance~\cite{STREAM} is 11 gigabytes per second of achieved
memory bandwidth. All tests are run on a single processor core.

\begin{figure}[ht]
\centering
{\includegraphics[width=0.45\textwidth]{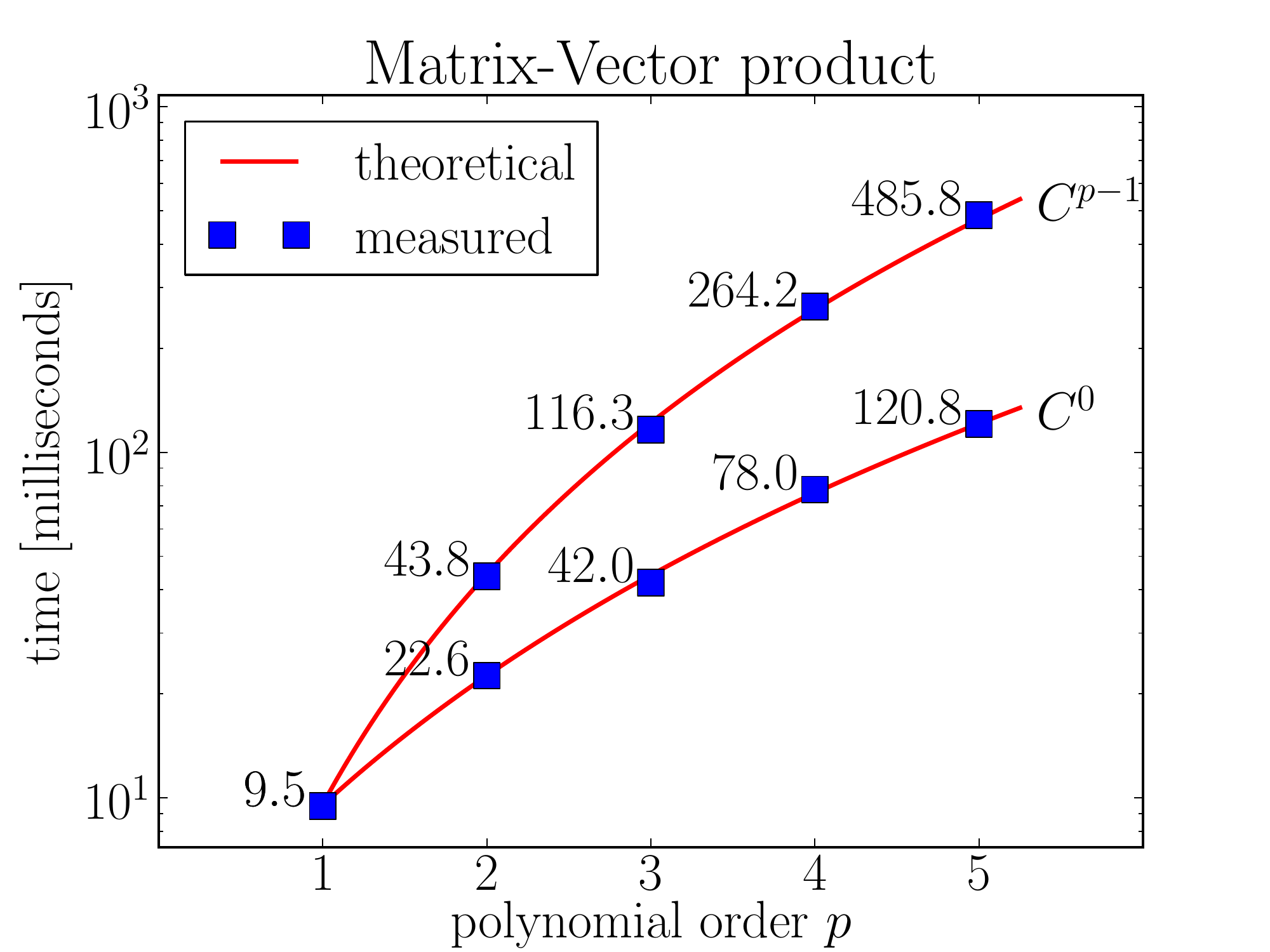}}
{\includegraphics[width=0.45\textwidth]{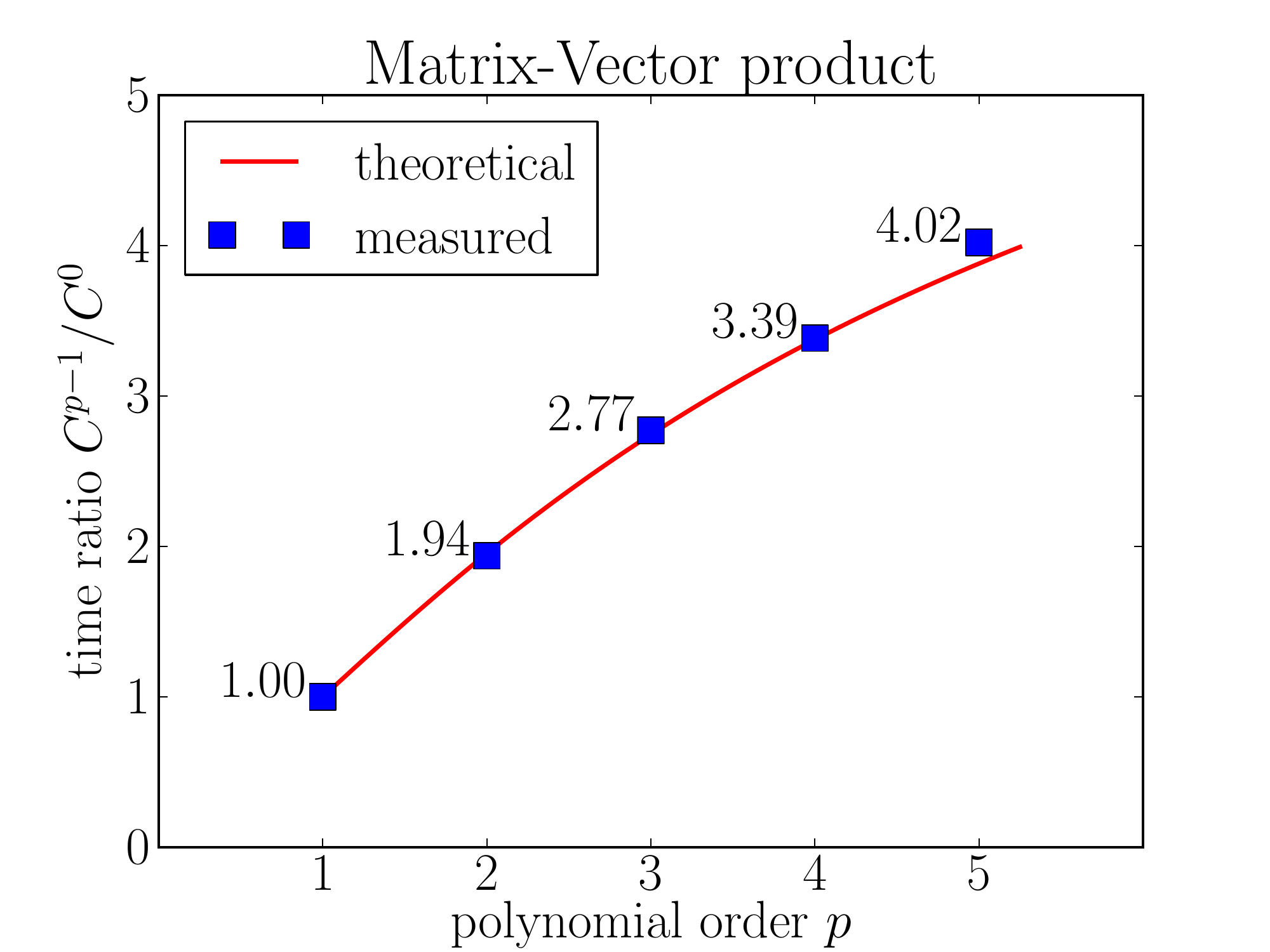}}
\caption{Matrix-Vector product}\label{f:SpMV}
\end{figure}

\begin{figure}[ht]
\centering
{\includegraphics[width=0.45\textwidth]{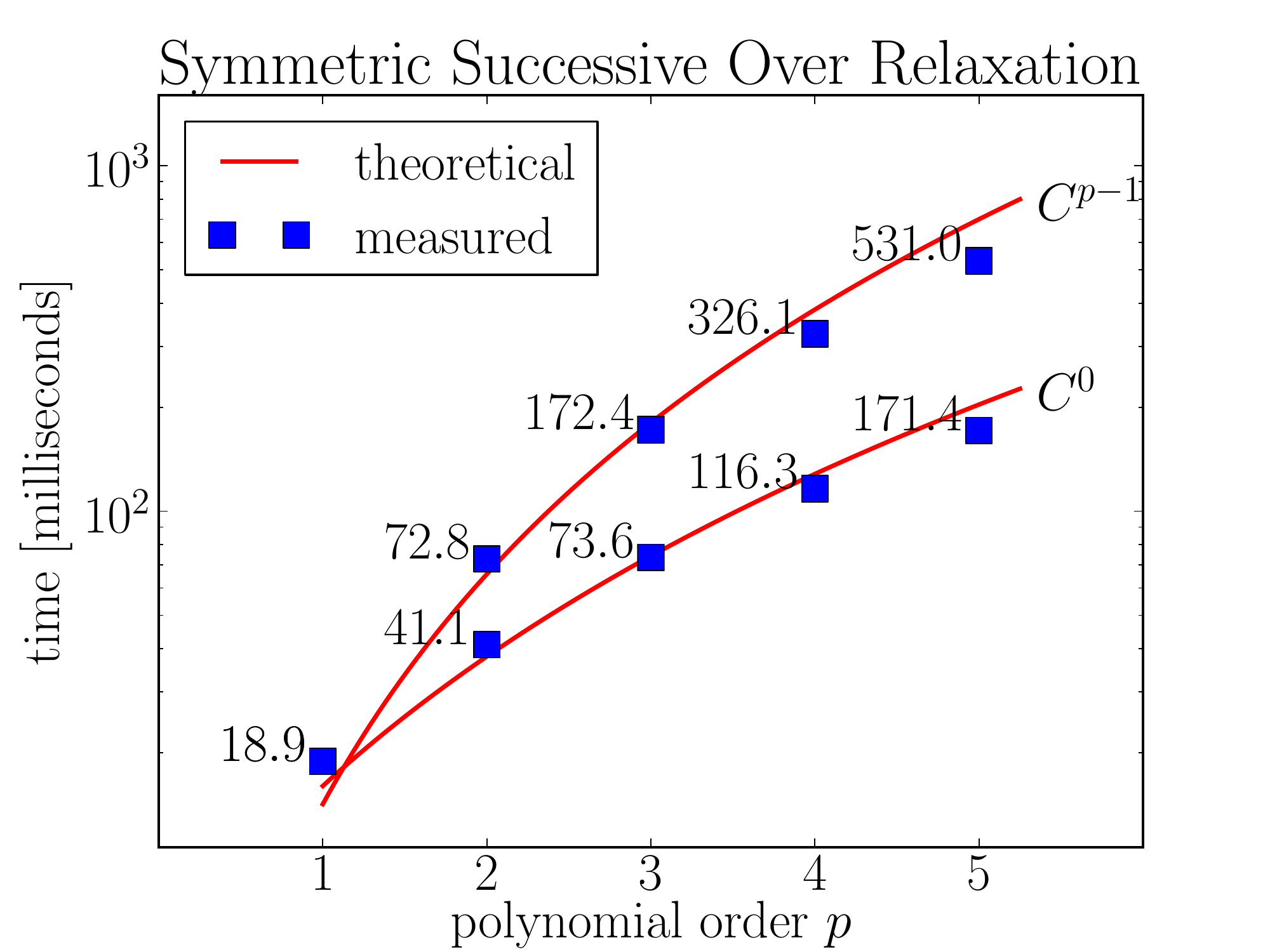}}
{\includegraphics[width=0.45\textwidth]{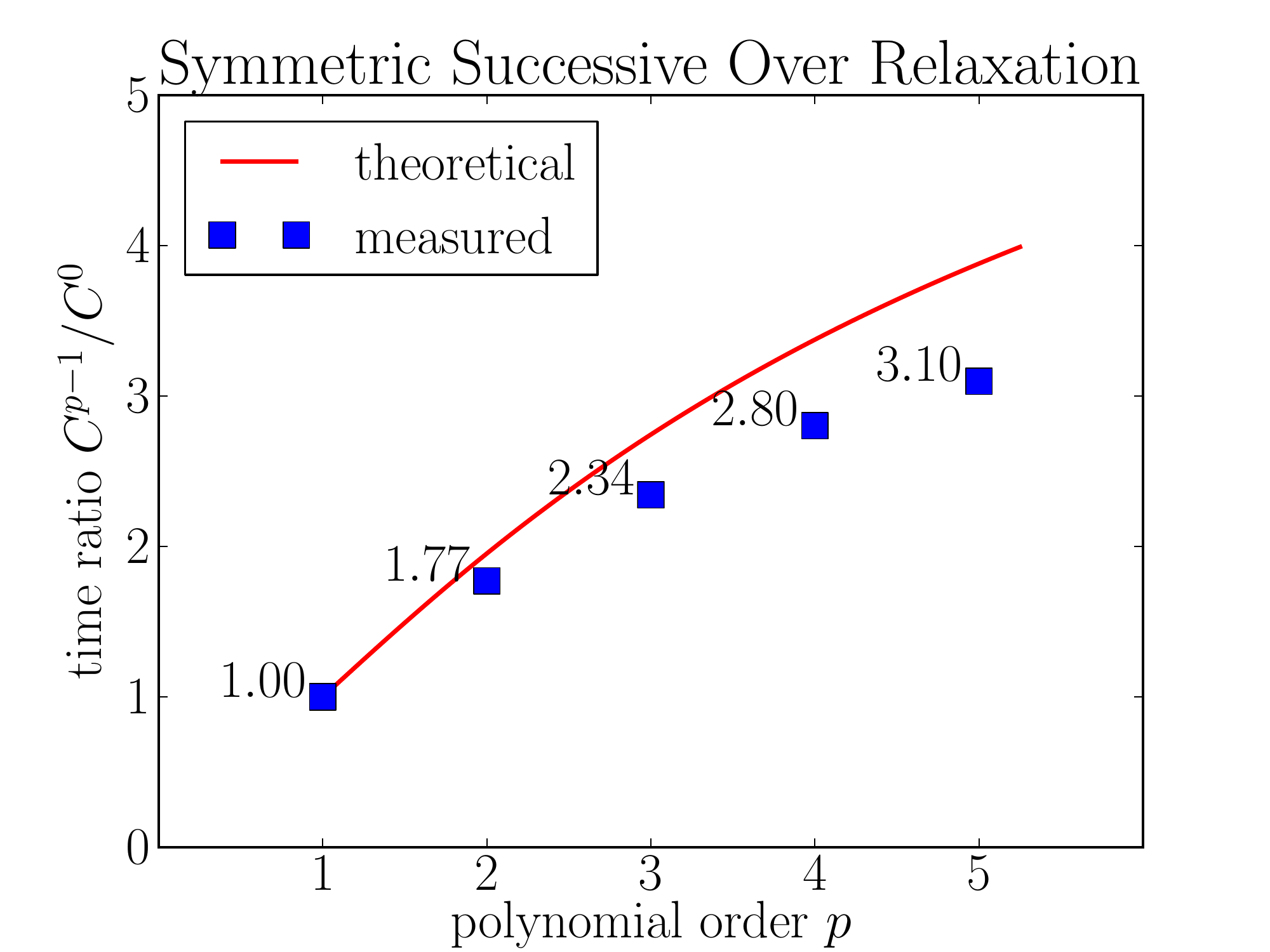}}
\caption{Symmetric Successive Over Relaxation}\label{f:SSOR}
\end{figure}

\begin{figure}[ht]
\centering
\subfloat[factorization]{%
{\includegraphics[width=0.45\textwidth]{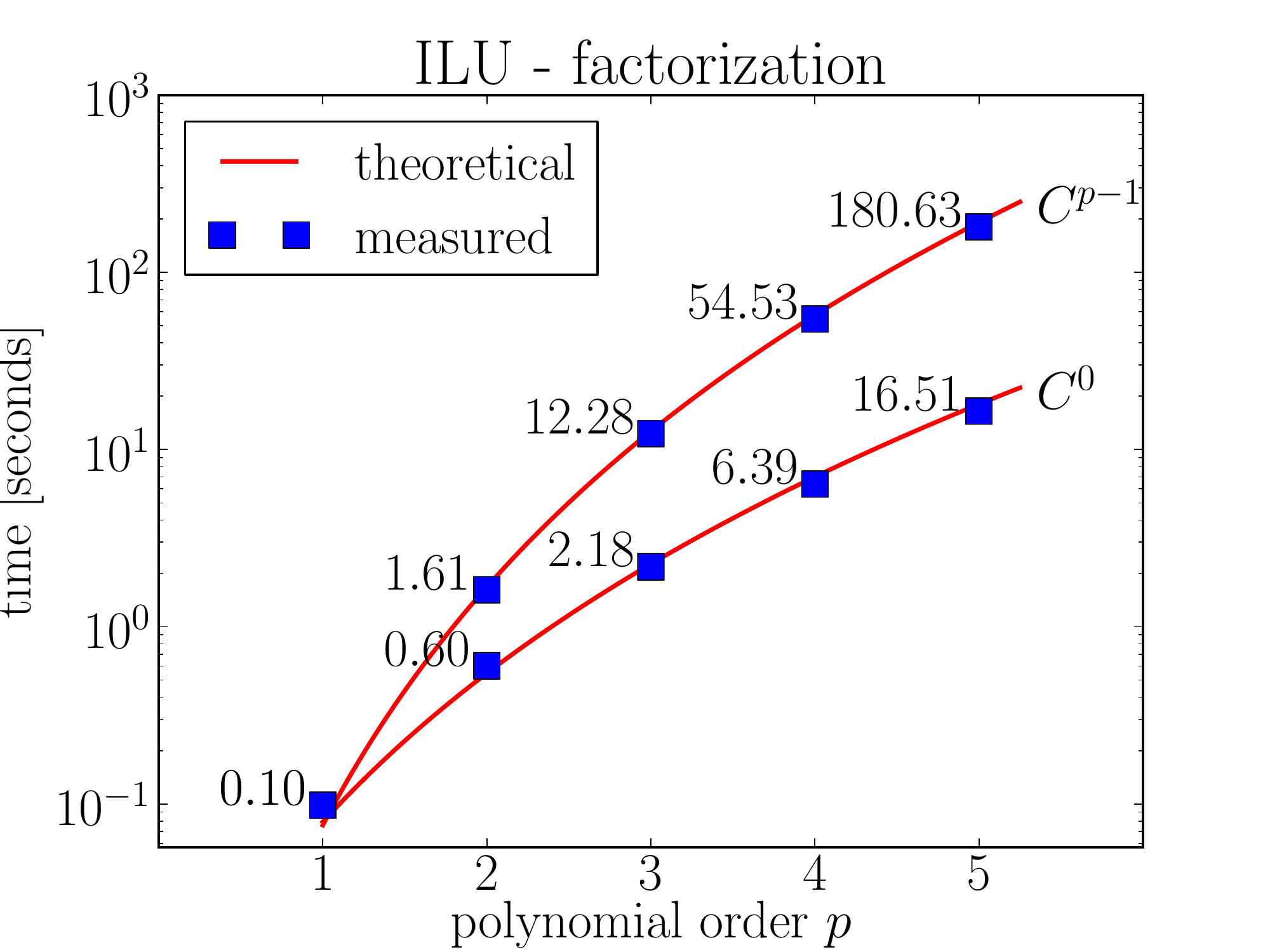}}
{\includegraphics[width=0.45\textwidth]{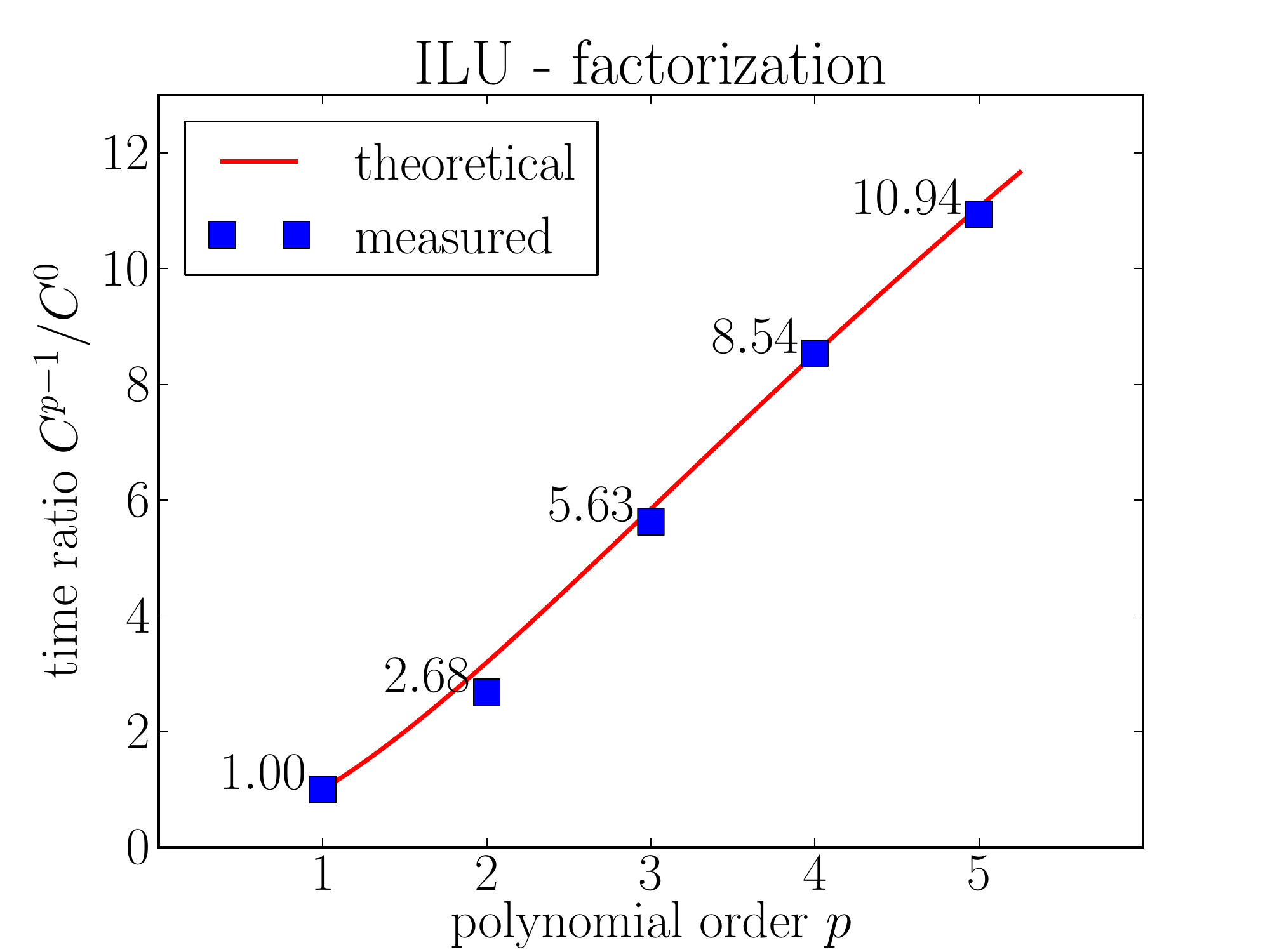}}
}\\
\subfloat[triangular solves]{%
{\includegraphics[width=0.45\textwidth]{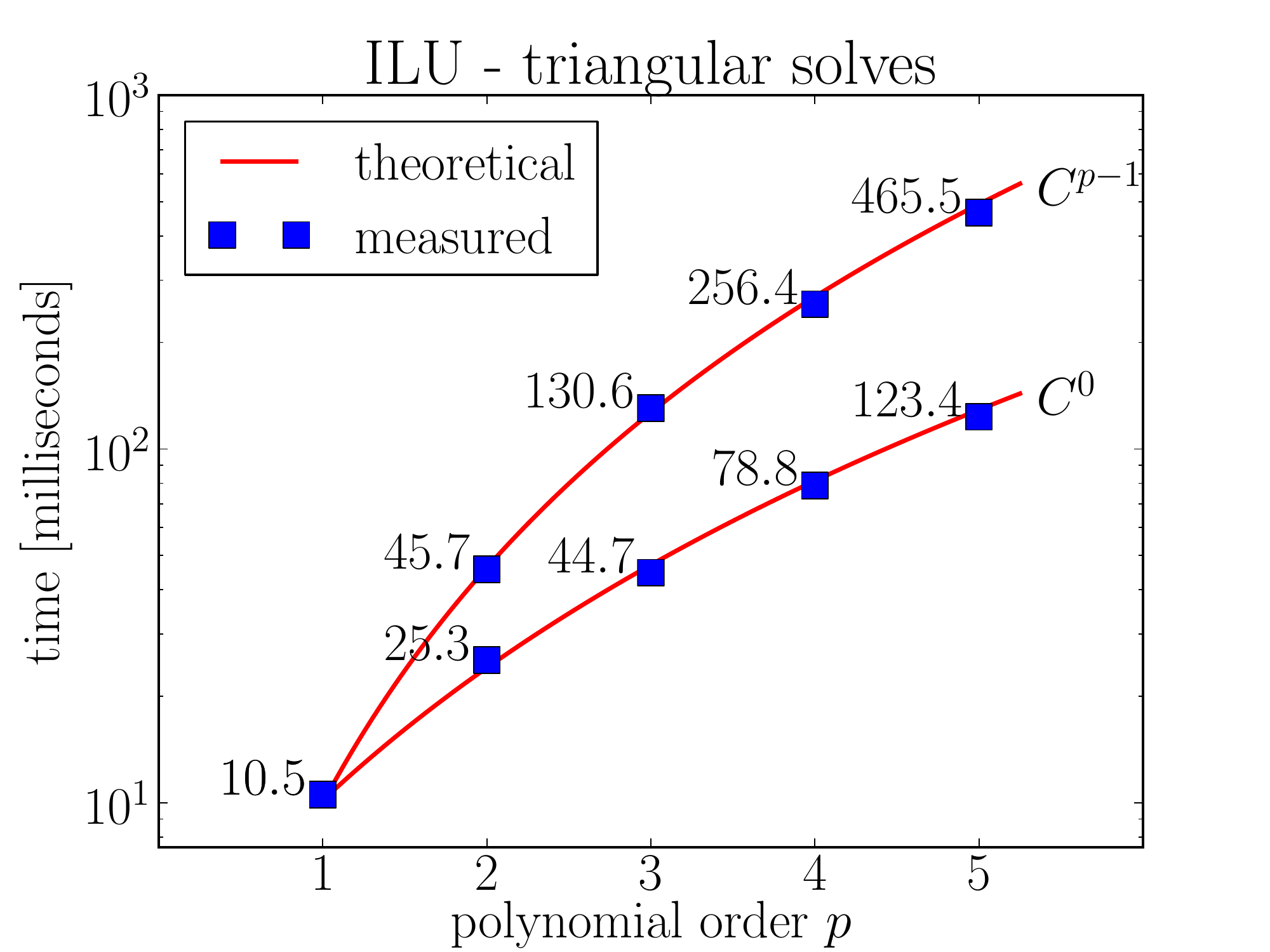}}
{\includegraphics[width=0.45\textwidth]{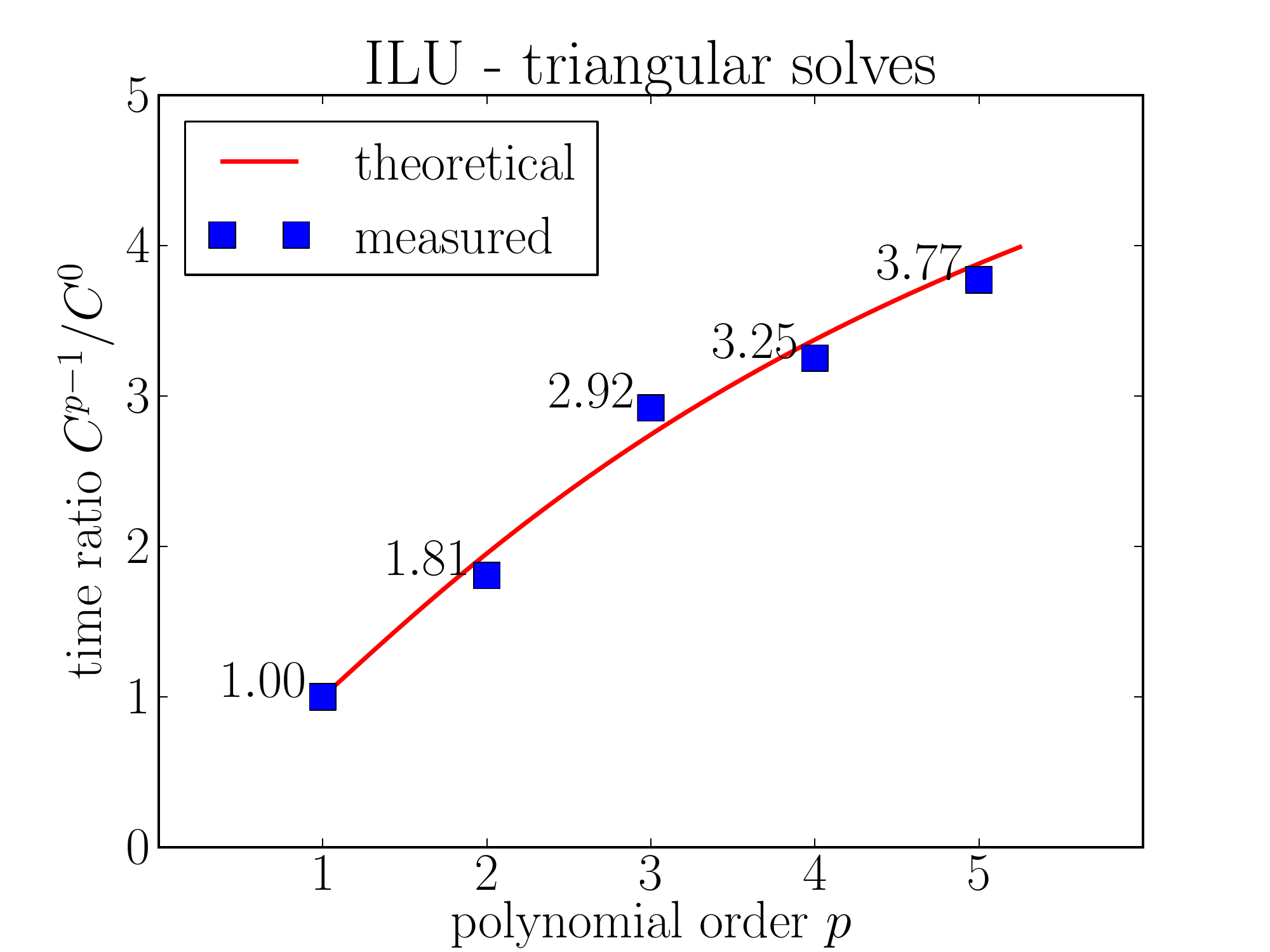}}
}
\caption{Incomplete LU}\label{f:ILU}
\end{figure}

Figures~\ref{f:SpMV}, ~\ref{f:SSOR}, and ~\ref{f:ILU} present the
results of our experiments. Plots on the left present estimated and
measured wall-clock times for $C^{0}$ and $C^{p-1}$ spaces. Square
markers correspond to measured time, solid lines correspond to our
theoretical FLOP count estimates scaled with the achieved mean FLOP
rates in order to relate FLOP count to time. Plots on the right
present the time ratio for $C^0$ and $C^{p-1}$ discretizations.

Overall, the measurements match closely the theoretical estimates
except for SSOR sweeps. In the case of SSOR, forward and backward
sweeps operate on the lower and upper triangular parts of the sparse
matrix. Furthermore, the backward sweep is performed in reversed row
ordering. The data access pattern is much more irregular than the one
in matrix-vector product and causes greater miss rates in the
processor cache. This situation leads to lower FLOP rates, hindering
performance. Textbook implementations of ILU factorization and
forward/backward triangular solves also suffer from this
issue. However, PETSc employs a different data layout to store the LU
factors and is able to achieve a FLOP rate comparable to the one for
matrix-vector products. See~\cite{sz2010} for a thorough analysis and
discussion about the importance of data layout in triangular solves.

\subsection{Iteration Counts}

The purpose of a preconditioner is to improve the spectrum of the
eigenvalues of the original operator. The number of iterations
required for convergence is tightly related to this
spectrum~\cite{Nachtigal1992}. While the setup and application cost of
the preconditioner is an important factor, a measure of how well
preconditioners work in terms of number of iterations is also
critical. In tables~\ref{t:itC0} and \ref{t:itCpm1}, we present
numerical results which test how well each preconditioner works in
terms of number of iterations.

We study the iteration counts as the function spaces vary in $h$ and
$p$. In this study we interpret $h$ as half of the basis support size
as opposed to the traditional interpretation of the element size,
denoted as $h_e$. This means that $h$ retains its original meaning in
the case of $C^0$ spaces, that is $h=h_e$. However, for $C^{p-1}$
spaces, $h=\slfrac{h_e(p+1)}{2}$. We argue this based on the
observation that under this definition, the condition number scales as
standard theory suggests ($h^{-2}$, \cite{Axelsson2003}) for the
Laplace problem. If one considers scaling of the condition number
based on $h_e$, the integration support, then the scaling artificially
appears to be better.
\begin{table}
\caption{Number of iterations required for convergence of CG using
  different preconditioners and $C^0$ B-spline spaces}\label{t:itC0}
\begin{center}
\begin{tabular}{clcccc}
\hline
   & & \multicolumn{4}{c}{Basis support size, $h$}\\
\cline{3-6}
$p$ & type & $\slfrac{1}{2}$ & $\slfrac{1}{4}$ & $\slfrac{1}{8}$ & $\slfrac{1}{16}$ \\
\hline
1 & Jacobi & 5   & 11  & 23  & 45 \\
2 & Jacobi & 19  & 31  & 42  & 61   \\
3 & Jacobi & 68  & 91  & 97  & 107   \\
4 & Jacobi & 216 & 355 & 406 & 424  \\
\rowcolor[gray]{.9}1 & SSOR & 5 & 7 & 14 & 24 \\
\rowcolor[gray]{.9}2 & SSOR & 11 & 14 & 17 & 27 \\
\rowcolor[gray]{.9}3 & SSOR & 27 & 30 & 31 & 37 \\
\rowcolor[gray]{.9}4 & SSOR & 65 & 75 & 77 & 80 \\
1 & ILU & 1 & 7 & 12 & 22 \\
2 & ILU & 7 & 9 & 15 & 26 \\
3 & ILU & 8 & 11 & 19 & 34 \\
4 & ILU & 10 & 14 & 24 & 42 \\
\rowcolor[gray]{.9}1 & EBE &  6  & 21  & 30  & 45 \\
\rowcolor[gray]{.9}2 & EBE & 17  & 34  & 45  & 65 \\
\rowcolor[gray]{.9}3 & EBE & 23  & 38  & 50  & 86 \\
\rowcolor[gray]{.9}4 & EBE & 27 & 47 & 63 & 111   \\
\hline
\end{tabular}
\end{center}
\end{table}
\begin{table}
\caption{Number of iterations required for convergence of CG using different preconditioners and $C^{p-1}$ B-spline spac
es. The BBB preconditioner is shown for $r=\slfrac{p}{2}$.}\label{t:itCpm1}
\begin{center}
\begin{tabular}{clcccc}
\hline
   & & \multicolumn{4}{c}{Basis support size, $h$}\\
\cline{3-6}
$p$ & type & $\slfrac{1}{2}$ & $\slfrac{1}{4}$ & $\slfrac{1}{8}$ & $\slfrac{1}{16}$ \\
\hline
1 & Jacobi & 5   & 11  & 23  & 45 \\
2 & Jacobi & 22  & 30  & 39  & 65 \\
3 & Jacobi & 70  & 99  & 100 & 105 \\
4 & Jacobi & 165 & 139 & 149 & 152 \\
\rowcolor[gray]{.9}1 & SSOR & 5 & 7 & 14 & 24 \\
\rowcolor[gray]{.9}2 & SSOR & 13 & 15 & 19 & 29 \\
\rowcolor[gray]{.9}3 & SSOR & 33 & 34 & 34 & 41 \\
\rowcolor[gray]{.9}4 & SSOR & 82 & 68 & 65 & 66 \\
1 & ILU & 1 & 7 & 12 & 22 \\
2 & ILU & 5 & 7 & 12 & 22 \\
3 & ILU & 6 & 7 & 12 & 20 \\
4 & ILU & 6 & 8 & 12 & 20 \\
\rowcolor[gray]{.9}1 & EBE &  6  & 21  & 30  & 45 \\
\rowcolor[gray]{.9}2 & EBE & 26  & 45  & 51  & 60  \\
\rowcolor[gray]{.9}3 & EBE & 50 & 74 & 81 &  89  \\
\rowcolor[gray]{.9}4 & EBE & 77 & 109 & 116 & 123 \\
1 & BBB & 5 & 11 & 23 & 45 \\
2 & BBB & 18 & 22 & 26 & 42 \\
3 & BBB & 30 & 35 & 38 & 54 \\
4 & BBB & 29 & 32 & 36 & 51 \\
\hline
\end{tabular}
\end{center}
\end{table}

In both spaces, the remarkable result is that the ILU preconditioner
outperforms other options in terms of number of iterations as well as
the $p$ scaling. Of greater interest is that in the case of $C^{p-1}$
spaces, its $p$ scalability is perfect. While no theory currently
exists to prove that the ILU preconditioner will lead to a converged
result, it is among the most economical preconditioners to setup and
apply to linear systems in the range of problems solved. 

Our interest in studying the cost of solving medium-size problems
using standard techniques is at the core of more complex
preconditioning approaches such as domain-decomposition and
multigrid. For example, despite the fact that ILU does not scale as
well in $h$, we can use it in a multigrid approach to remove $h$
dependence. For example, in table~\ref{t:itMG} we show iteration
counts for a two grid solver we constructed. On the fine grid, we use
CG with ILU and a direct solver on the coarse level. The coarse level
is a factor of $2^3$ unrefined in $h$ from the fine level. We see that
as in standard $C^0$ spaces, a multigrid approach is able to remove
$h$ dependence from linear systems discretized using $C^{p-1}$ spaces.

\begin{table}
\caption{Number of iterations required for convergence with a two grid
  solver}\label{t:itMG}
\begin{center}
\begin{tabular}{clcccc}
\hline
   & & \multicolumn{4}{c}{Basis support size, $h$}\\
\cline{3-6}
$p$ & space & $\slfrac{1}{2}$ & $\slfrac{1}{4}$ & $\slfrac{1}{8}$ & $\slfrac{1}{16}$ \\
\hline
1 & $C^0$ & 4  & 16  & 19  & 19 \\
2 & $C^0$ & 15  & 18  & 20  & 21 \\
3 & $C^0$ & 18  & 20  & 21  & 22 \\
4 & $C^0$ & 20  & 24  & 25  & 26 \\
\rowcolor[gray]{.9}1 & $C^0$ & 4 & 16 & 19 & 19 \\
\rowcolor[gray]{.9}2 & $C^1$ & 14 & 15 & 17 & 19 \\
\rowcolor[gray]{.9}3 & $C^2$ & 14 & 15 & 18 & 19 \\
\rowcolor[gray]{.9}4 & $C^3$ & 14 & 15 & 18 & 19 \\
\hline
\end{tabular}
\end{center}
\end{table}

The iteration counts for the BBB preconditioner shown in
table~\ref{t:itCpm1} are for an overlap parameter
$r=\slfrac{p}{2}$. We have selected this size in an attempt to balance
the cost of applying the preconditioner and its convergence. In
table~\ref{t:bbb}, we show convergence results for more choices of the
overlap parameter $r$. Note that when $r=0$ the preconditioner is
diagonal Jacobi. As the overlap increases, we see an improvement in
the number of iterations. However, when the overlap is at its maximum,
$r=p$, the setup and application of the preconditioner is
prohibitively expensive. We suggest that the choice $r=\slfrac{p}{2}$
leads to a good compromise between fast convergence and moderate
application cost.

\begin{table}[ht]
  \caption{Convergence results for the BBB preconditioner where the
    basis support size a constant $h=0.25$ and the overlap $r$
    varies. Underlined entries represent a preconditioner with
    approximately the same number of nonzero entries as the original
    system matrix.}\label{t:bbb}
\begin{center}
\begin{tabular}{cccccc}
\hline
   & \multicolumn{5}{c}{Basis overlap, $r$}\\
\cline{2-6}
$p$ & 0 & 1 & 2 & 3 & 4\\
\hline
1  & \underline{11} & 16 &  & &  \\
2  & 30 & \underline{22} & 21 & &  \\
3  & 99 & \underline{35} & 25 & 24 &  \\
4  & 139 & 56 & \underline{32} & 28 & 27 \\
5  & 345 & 76 & \underline{48} & 35 & 30 \\
6  & 353 & 112 & 68 & \underline{48} & 35 \\
7  & 610 & 167 & 104 & \underline{64} & 43 \\
8  & 737 & 201 & 130 & 92 & \underline{65} \\
\hline
\end{tabular}
\end{center}
\end{table}

\subsection{Timing Results}

While the number of iterations required for convergence is a useful
measure to study, it is not sufficient to understand which
preconditioning options are better from a practical point of
view. Frequently, a practitioner must experiment with different options
on a meaningful range of problem sizes. For a preconditioner to be
effective, the cost of setting up and applying must be weighted
against its capability to improve the spectrum of eigenvalues,
effectively reducing the total number of iterations. 

\begin{figure}[ht]
\centering
\subfloat[$p=2$]{\includegraphics[width=0.45\textwidth]{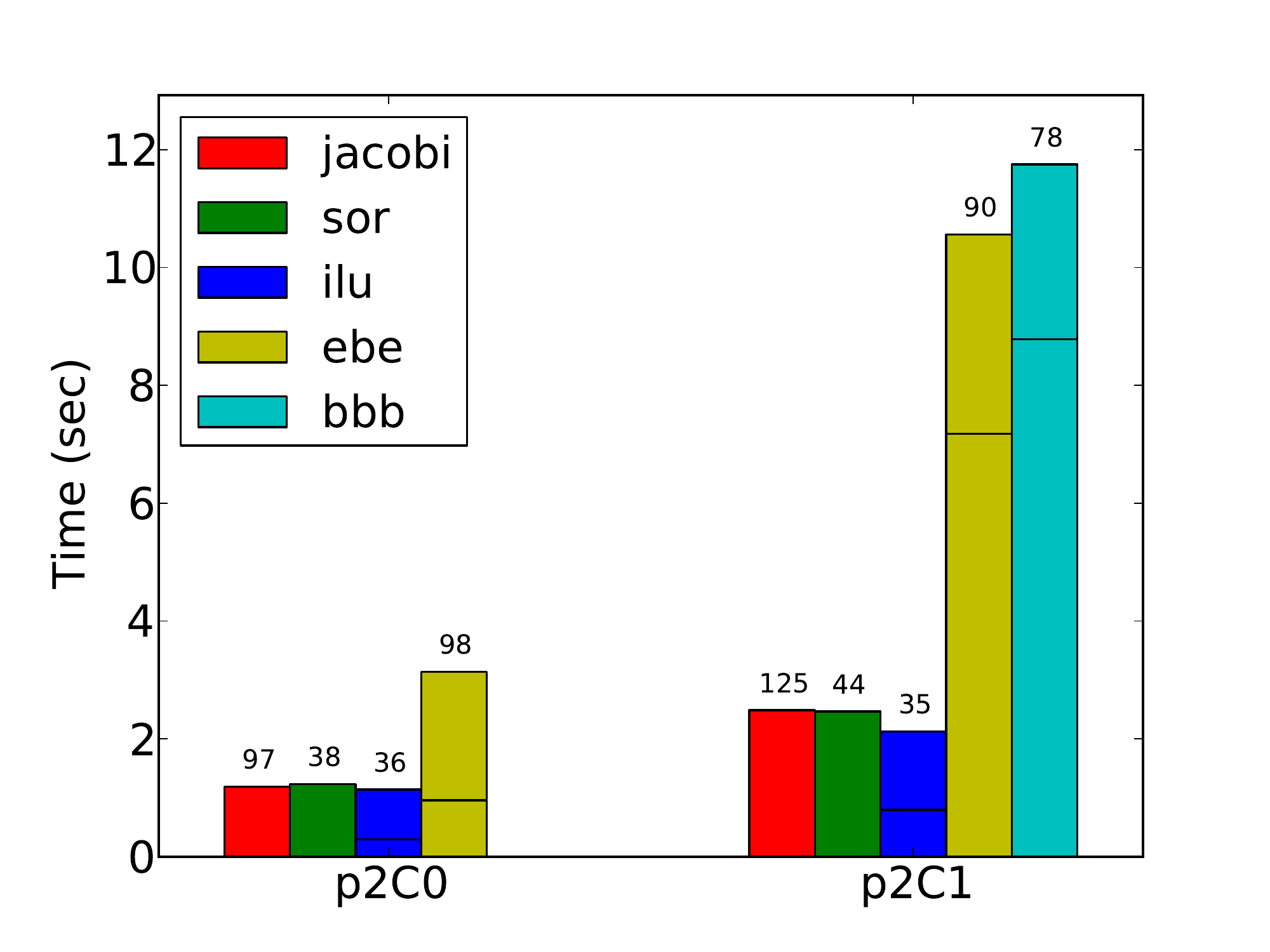}}
\subfloat[$p=3$]{\includegraphics[width=0.45\textwidth]{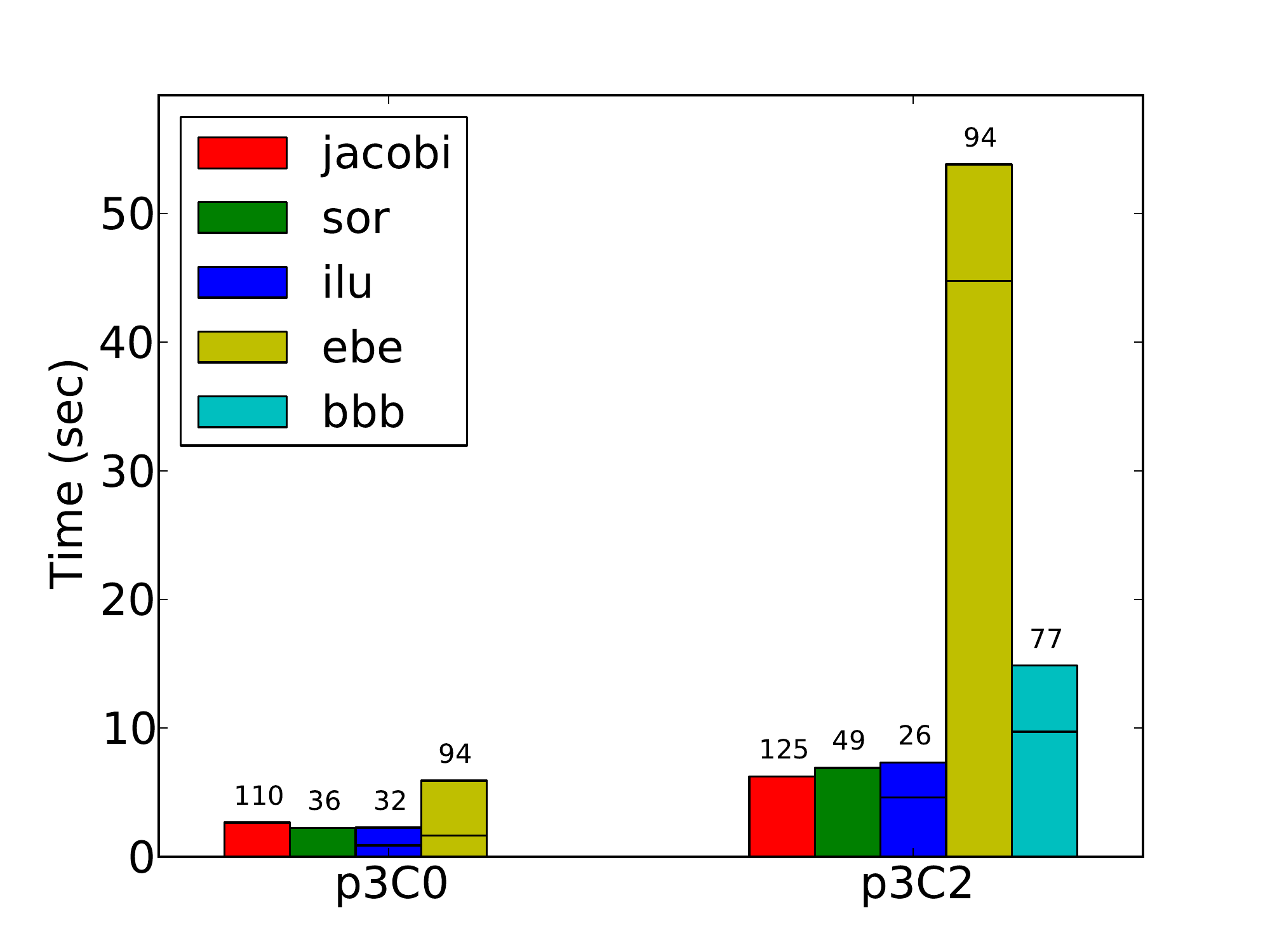}}\\
\subfloat[$p=4$]{\includegraphics[width=0.45\textwidth]{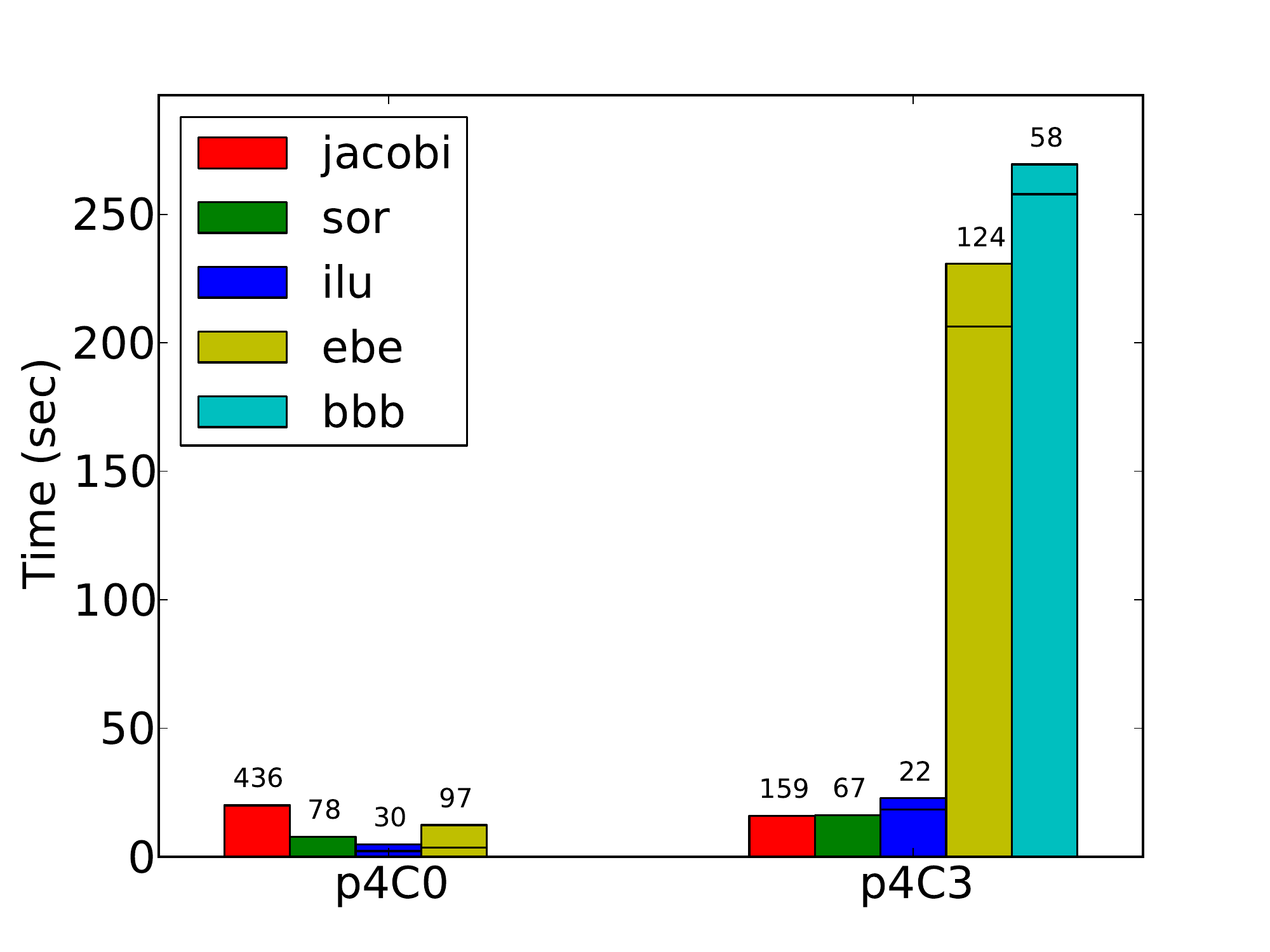}}
\subfloat[$p=5$]{\includegraphics[width=0.45\textwidth]{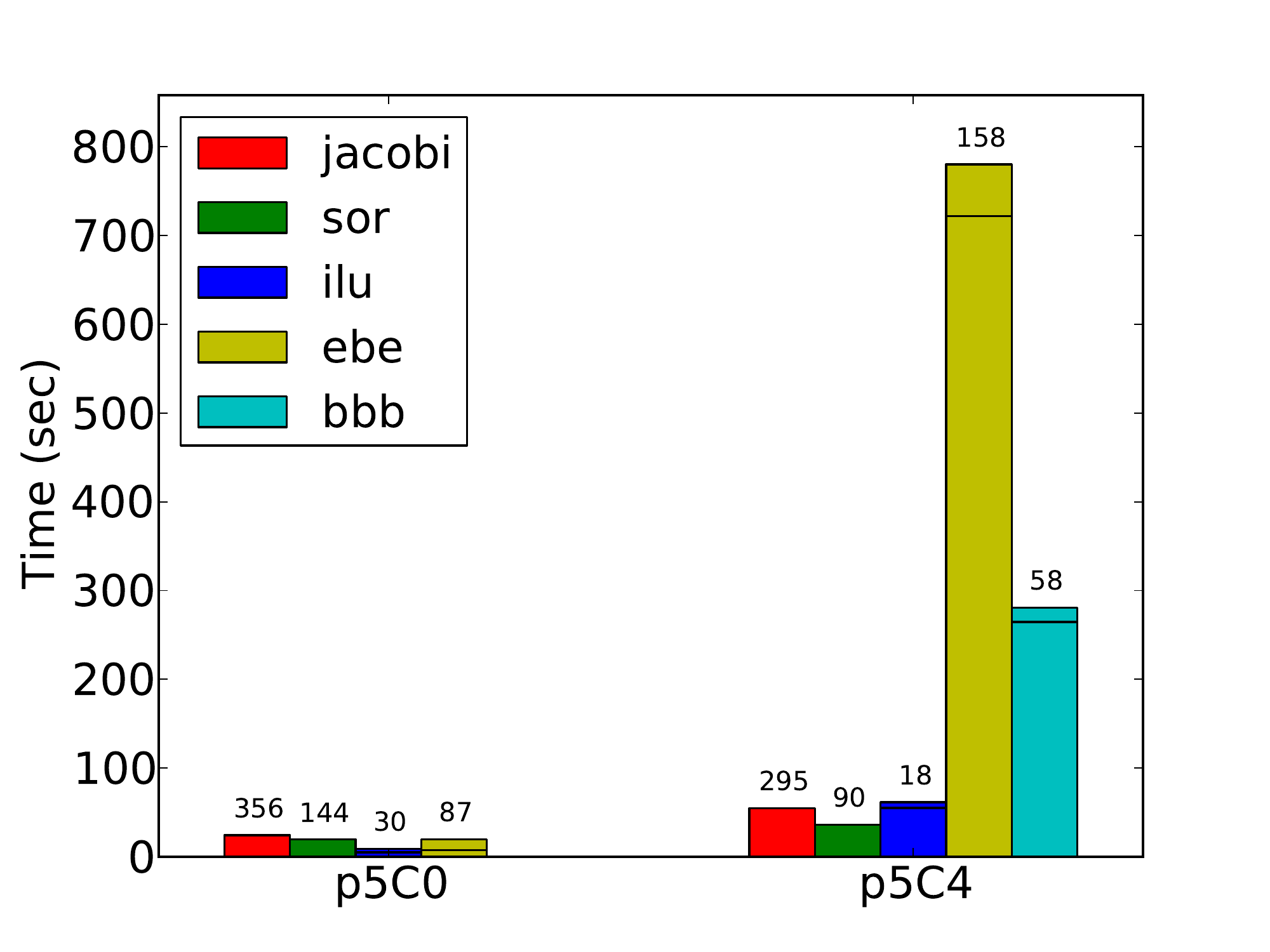}} 
\caption{Solution times of various preconditioning options for $C^0$ and $C^{p-1}$ spaces consisting of $N=10^5$ degrees of freedom.}\label{f:100k}
\end{figure}

We first present some timing results for linear systems consisting of
$10^5$ degrees of freedom. In figure~\ref{f:100k} we display bar
graphs representing the total solution time required for convergence
for different preconditioning options and varying polynomial order
$p$. In each plot we display $C^0$ spaces on the left and $C^{p-1}$ on
the right. Furthermore, each bar is divided into two parts. The bottom
part represents the setup time required for each preconditioner and
the top the remaining solve time. We also include the number of
iterations required for convergence on the top of each bar.

Most striking is the time required by the EBE and BBB preconditioners
for $C^{p-1}$ spaces. As predicted in the theoretical estimates for
setup cost, these preconditioners are considerably more expensive to
setup, and additionally they are not able to significantly reduce the
total number of iterations. This is an effect that continues to grow
as we increase $N$, the size of the problem.

We would like to make an additional comment on the ILU
preconditioner. From the analysis of setup costs, it is clear that EBE
and BBB preconditioners are $p^3$ times more expensive to setup than
ILU. However, the EBE, BBB, and ILU preconditioners are all based on
the Gaussian elimination process. The EBE and BBB preconditioners are
built by taking into account the interaction between a compact
neighborhood of DOF (defined by elements in the case of EBE and by the
overlap $r$ in the case of BBB). Due to its algorithmic structure, we
intuitively argue that ILU is able to capture the interactions between
different DOF in a more global manner, and therefore has a better
effect on the convergence of the problem. 

\begin{figure}[ht]
\centering
\subfloat[$p=2$]{\includegraphics[width=0.45\textwidth]{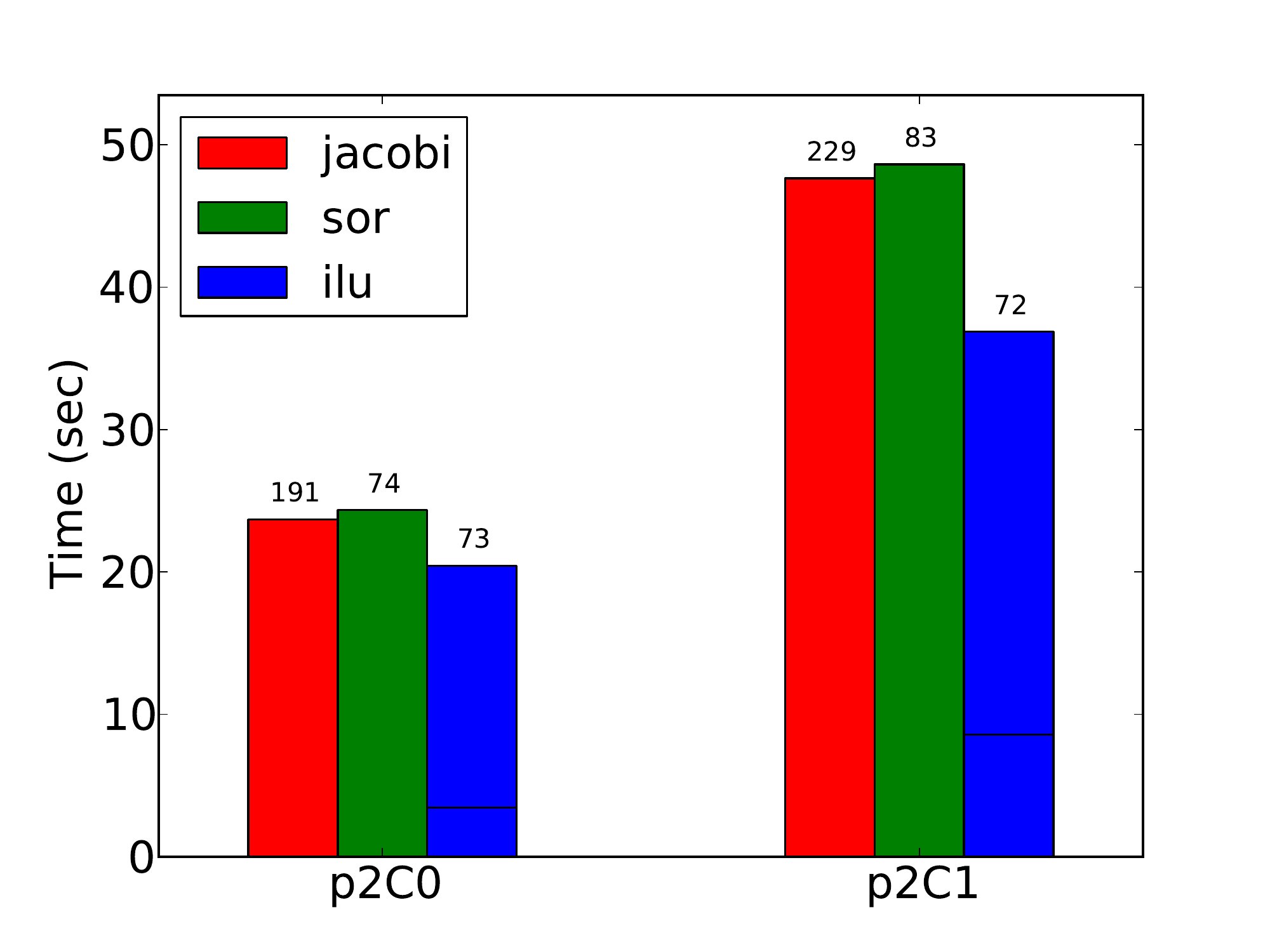}}
\subfloat[$p=3$]{\includegraphics[width=0.45\textwidth]{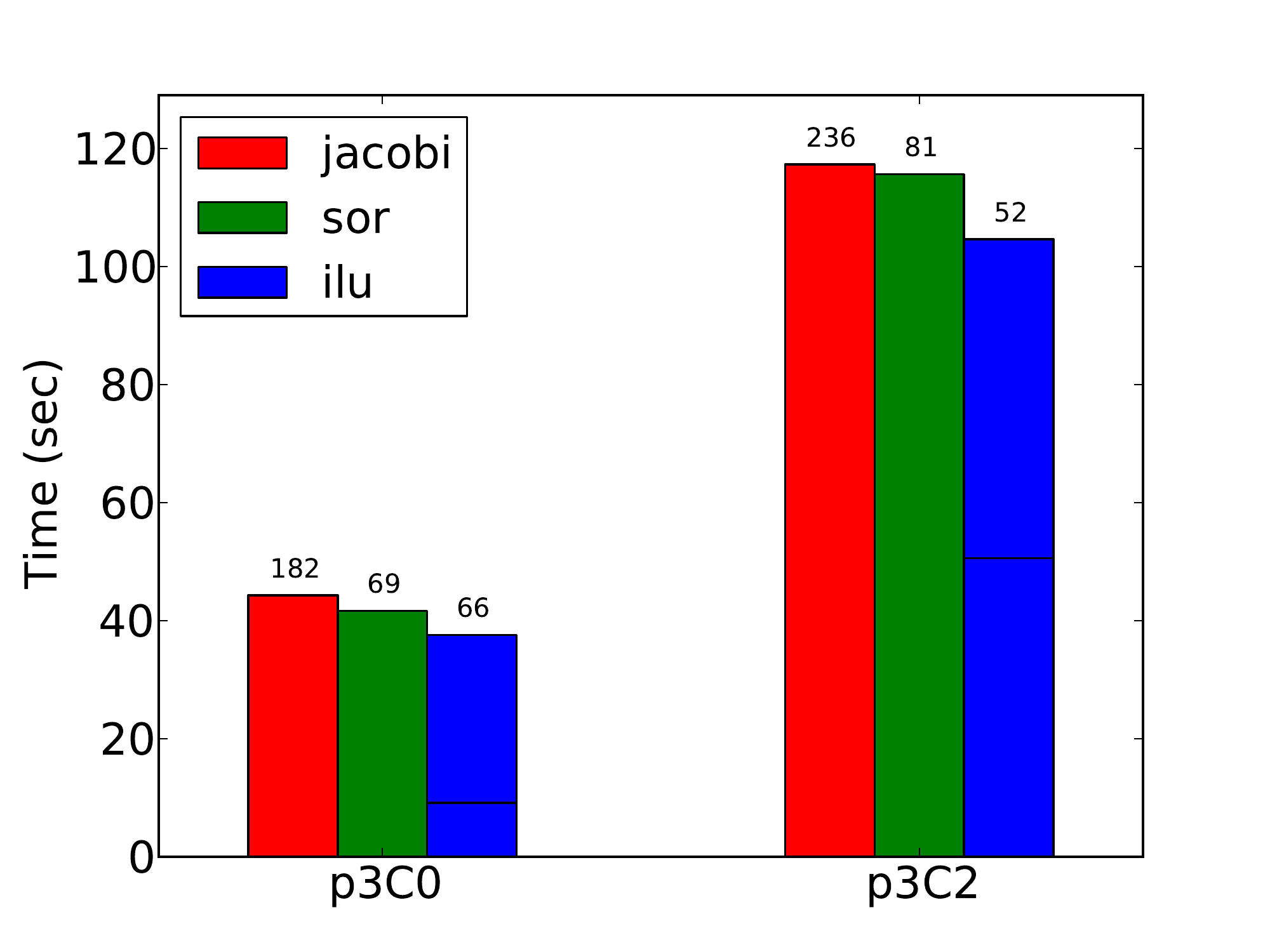}}
\caption{Solution times of various preconditioning options for $C^0$ and $C^{p-1}$ spaces consisting of $N=10^6$ degrees of freedom.}\label{f:1000k}
\end{figure}

In figure~\ref{f:1000k} we remove the EBE and BBB preconditioning to
better highlight the differences between the remaining choices. Also,
we have increased the number of DOF to $10^6$. We see that in this
case, the remaining preconditioners (Jacobi, SSOR, and ILU) all
perform similarly in terms of time. We also see that SSOR and ILU
require a similar number of iterations.
We also note that the $C^{p-1}$ spaces are two times as expensive than
the $C^0$ spaces for $p=2$ and three times as expensive for $p=3$, as
predicted by our theoretical estimates.

\begin{figure}[ht]
\centering
\subfloat[$C^{p-1}$ versus $C^0$]{\includegraphics[width=0.8\textwidth]{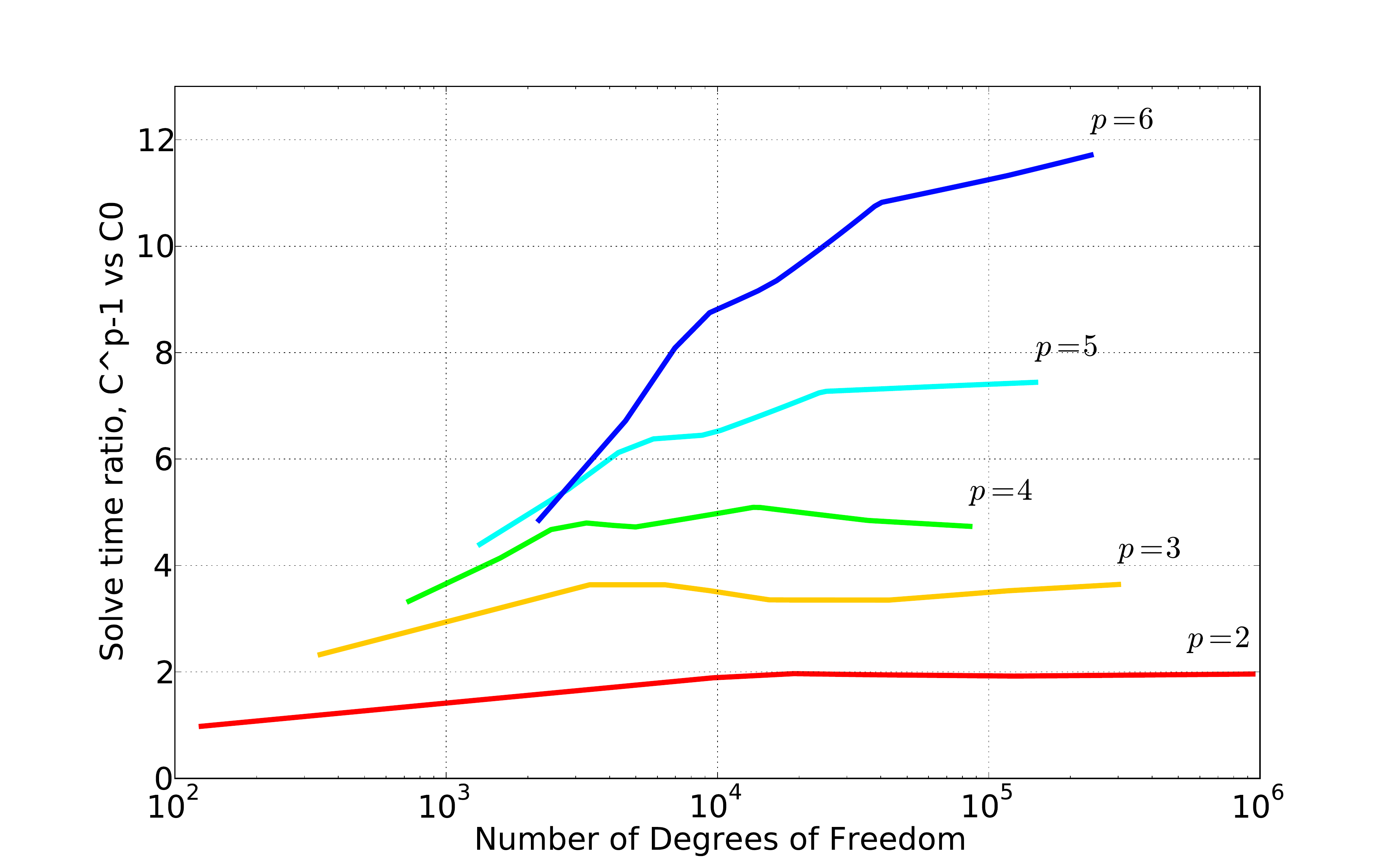}\label{f:ilu}}\\
\subfloat[$C^{p-1}$ versus $C^0$ with static condensation]{\includegraphics[width=0.8\textwidth]{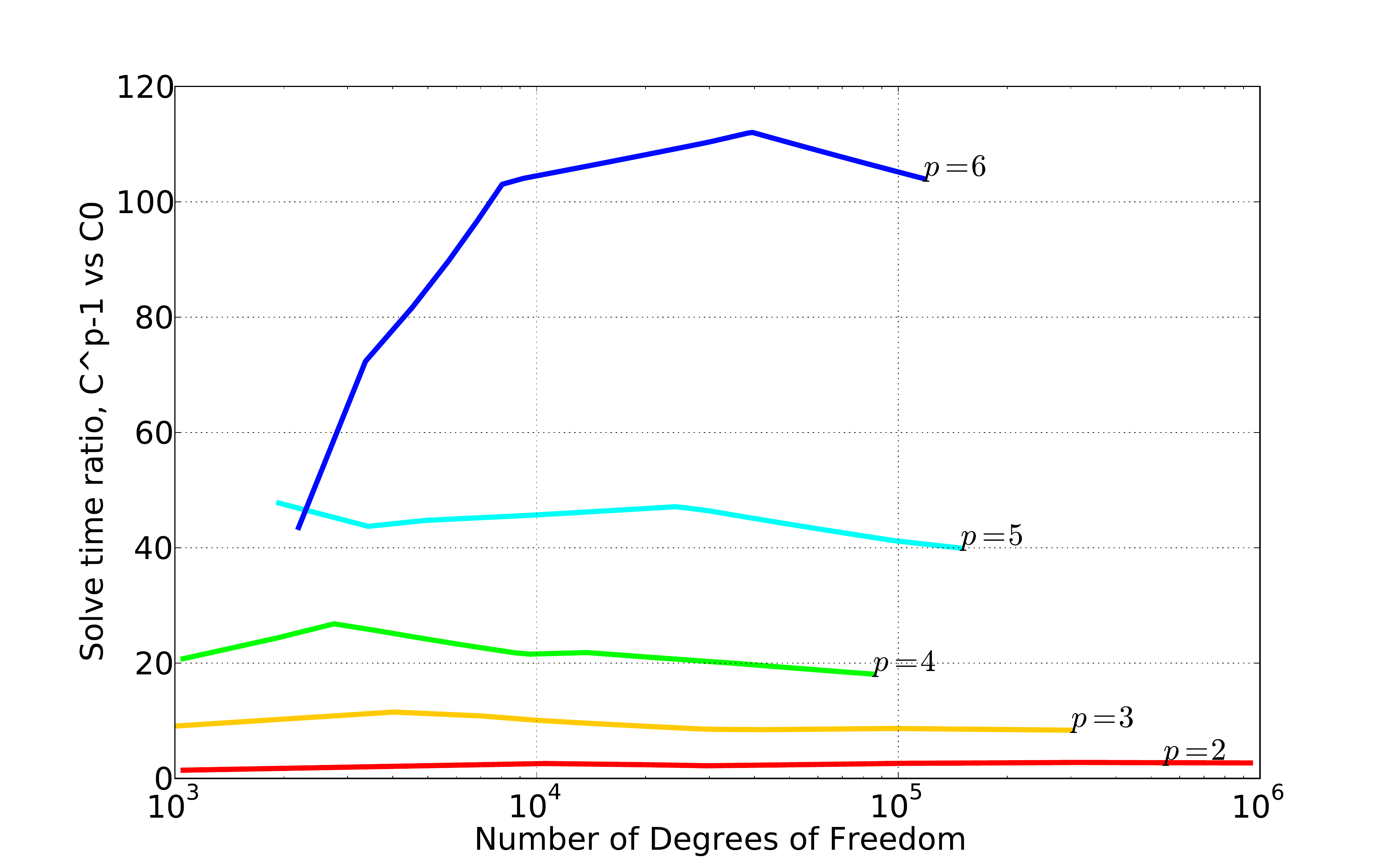}\label{f:ilusc}}
\caption{Solve time ratios for $C^{p-1}$ and $C^0$ spaces using CG preconditioned with ILU. Solve time is inclusive of the setup of the preconditioner but excludes assembly time. In the case of static condensation, the solve time includes extra assembly time required to compute the Schur complements.}
\end{figure}

We extend the study to include a range of degrees of freedom in
figure~\ref{f:ilu}, however we only consider the ILU
preconditioner. We show the ratio of solve times for $C^{p-1}$ to
$C^0$ for a range of problem sizes $N$ and polynomial orders $p$. We
highlight that for this example, the cost of the use of the higher
continuous basis is anywhere from $p$ to $2p$ times more
expensive. However, in the theory section, we showed that the solve
cost can be greatly reduced when using static condensation and solving
the skeleton problem. In figure~\ref{f:ilusc} we have repeated the
numerical experiment, this time using static condensation on the $C^0$
spaces. Note that the additional assembly time incurred due to static
condensation operations have been included into the solve cost. In
this case, the $C^{p-1}$ spaces are $\mathcal{O}(p^2)$ times more
expensive to solve, as predicted by the theoretical cost of the
matrix-vector multiplication estimates.

%% file: conc.tex
We have presented a study on the additional cost incurred in the
iterative solver due to the use of a more continuous basis in a
Galerkin weak form. We have presented theoretical estimates for
computational costs of matrix-vector multiplication as well as
preconditioner setup and application for a variety of preconditioning
techniques. We present numerical results for the Laplace problem to
establish a baseline understanding of how continuity affects the
solver. 

We conclude that the matrix-vector product is at most eight times more
expensive for the $C^{p-1}$ spaces. However, when using high $p$ with
static condensation, this factor increases to $\slfrac{8p^2}{33}$. We
expect that the improved approximability per DOF of the $C^{p-1}$
spaces may be better realized when using the iterative solver,
particularly for low $p$. This is, however, strongly tied to the
performance of the selected preconditioner applied to the equation of
interest.

We observe that for moderate $N$ and the Laplace problem, that the EBE
and BBB preconditioners are prohibitively expensive options for
$C^{p-1}$ spaces. This is because there is an element per basis
function which results in far greater setup costs than in $C^0$
spaces. For these options to be effective for a particular problem,
they should lead to a considerably large decrease in iterations
compared to other options. This is not the case for the simple Laplace
problem and, intuitively, we do not expect this to be the case for
more complex applications.

The ILU preconditioner, while lacking a theoretical ground for
convergence, performs quite well in terms of iterations and
computational time. Remarkably, for $C^{p-1}$ spaces we observe
perfect almost $p$-scaling of the preconditioned operator condition
number up to $10^6$ degrees of freedom.